ORIGINAL ARTICLE

# Practicable robust stochastic optimization under divergence measures with an application to equitable humanitarian response planning

**Aakil M. Caunhye**[1] · **Douglas Alem**[1]

Received: 20 July 2022 / Accepted: 18 April 2023 / Published online: 21 May 2023

**Abstract**

We seek to provide practicable approximations of the two-stage robust stochastic optimization model when its ambiguity set is constructed with an *f*-divergence radius. These models are known to be numerically challenging to various degrees, depending on the choice of the *f*-divergence function. The numerical challenges are even more pronounced under mixed-integer first-stage decisions. In this paper, we propose novel divergence functions that produce practicable robust counterparts, while maintaining versatility in modeling diverse ambiguity aversions. Our functions yield robust counterparts that have comparable numerical difficulties to their nominal problems. We also propose ways to use our divergences to mimic existing *f*-divergences without affecting the practicability. We implement our models in a realistic location-allocation model for humanitarian operations in Brazil. Our humanitarian model optimizes an effectiveness-equity trade-off, defined with a new utility function and a Gini mean difference coefficient. With the case study, we showcase (1) the significant improvement in practicability of the robust stochastic optimization counterparts with our proposed divergence functions compared to existing *f*-divergences, (2) the greater equity of humanitarian response that the objective function enforces and (3) the greater robustness to variations in probability estimations of the resulting plans when ambiguity is considered.

**Keywords** Stochastic programming · Ambiguity · Robust stochastic optimization · Moreau-Yosida regularization · *f*-divergence · Equitable humanitarian logistics

✉ Aakil M. Caunhye
Aakil.Caunhye@ed.ac.uk

Douglas Alem
Douglas.Alem@ed.ac.uk

[1] Business School, The University of Edinburgh, 29 Buccleuch Place, Edinburgh EH8 9JU, UK

## 1 Introduction

Scenario-based stochastic programming is a modeling paradigm in which a decision maker implements a set of proactive here-and-now decisions before observing the realization of uncertain parameters, which arise in the form of probabilistically distributed scenarios, and scenario-dependent recourse decisions. In classical stochastic programming, the objective is to minimize the sum of the total cost incurred from here-and-now decisions and the expectation, over a *known* distribution, of the total cost incurred from recourse actions Shapiro et al. (2009). In reality, the probability distribution of scenarios is often evaluated based on unreliable or ambiguous information/data. This, combined with the underlying optimization routine, leads to highly biased perceived value-addition of “optimal" decisions and poor out-of-sample performances, a phenomenon conventionally referred to as the *optimizer’s curse* (Smith and Winkler 2006). Robust stochastic optimization is an alternative modeling paradigm where, instead of the expected recourse cost over a single scenario distribution, the worst-case expected recourse cost over a family of distributions is minimized. This family of distributions is modeled within what is known as an ambiguity set. The concept of ambiguity originates from seminal decision theory works (Knight 2012) where *ambiguity*, the uncertainty in the probability distribution of a parameter, is explicitly distinguished from *risk*, the uncertainty in the realization of a parameter from a known probability distribution.

The focal point of robust stochastic optimization is the ambiguity set, which is the vehicle through which attitudes towards ambiguity are expressed. There are two general ways to construct ambiguity sets: 1) using descriptive statistics and 2) using statistical distances. Construction using descriptive statistics is achieved by considering distributions for which certain statistics on primitive uncertainties satisfy pre-specified bounds. An ambiguity set built that way typically contains supports and first-order moments, together with more complex statistics such as covariances (Delage and Ye 2010; Goh and Sim 2010), confidence sets (Wiesemann et al. 2014), mean absolute deviations (Postek et al. 2018; Roos and den Hertog 2021), expectations of second-order-conic-representable functions (Bertsimas et al. 2018; Zhi 2017) and general second-order moments (Zymler et al. 2013). The resulting optimization model is generally a second-order conic or a semidefinite program. An ambiguity set constructed using statistical distances considers distributions that are within a certain chosen statistical distance from a nominal distribution. Among the statistical distances used in the literature are the Wasserstein metric (Esfahani and Kuhn 2018; Gao and Kleywegt 2022; Hanasusanto and Kuhn 2018; Zhao and Guan 2018; Luo and Mehrotra 2017; Saif and Delage 2021), the Prohorov metric (Erdoğan and Iyengar 2006) and general *f*-divergences (Ben-Tal et al. 2013; Love and Bayraksan 2014; Bayraksan and Love 2015; Jiang and Guan 2016; Liu et al. 2021), which include special cases such as the Kullback–Leibler function (Hu and Hong 2013) and the variation distance (Sun and Xu 2015).

Robust stochastic optimization specifically applies to decision-making problems with event-wise ambiguity, which means that the uncertain parameters have the additional characteristic of being associated with a discrete set of events. In several decision-making problems, uncertainty is inherently event-based. For instance, disasters, such as earthquakes and tsunamis, are event-based occurrences since historical data come from discrete incidents. In such cases, modeling ambiguity without event-wise considerations is fundamentally flawed because it avoids the correlations that exist between data from the same event. In addition, omitting event-wise associations in estimating descriptive statistics or statistical distances leads to the underlying mis-characterization that all events follow the same stochastic process.

In robust stochastic optimization, ambiguity is modeled via uncertainties in scenario/event probabilities. In the literature, ambiguity sets for robust stochastic optimization are most famously constructed by using an *f*-divergence as statistical distance measurement between scenario probabilities, and allowing probability variations such that their resulting total *f*-divergence does not exceed a pre-specified radius. From a decision-making standpoint, the choice of which *f*-divergence to employ is influenced by 1) the decision maker's attitude towards ambiguity and 2) the limit behaviours of *f*-divergence functions. A conservative decision maker, for instance, will prefer a pointwise smaller *f*-divergence function since it allows larger probability deviations within the same radius. The behaviours of *f*-divergence functions at 0 and $\infty$ determine whether scenarios can be "suppressed" (made impossible) or "popped" (made possible) (Love and Bayraksan 2014). While *f*-divergence is the more widely used distance measure for scenario-based robust stochastic optimization, the Wasserstein distance is favoured on situations where data is not scenario-based. Several recent works have proposed scenario-based robust stochastic optimization approaches under the Wasserstein distance (Gao and Kleywegt 2022; Xie 2020), but these appear to not be exactly comparable with *f*-divergence functions. *f*-divergence balls include only distributions that are absolutely continuous with respect to the nominal distribution. They thus do not accommodate distributions whose supports are bigger than the nominal distribution. This is usually used to motivate the use of the Wasserstein distance. In our case, however, this property of *f*-divergences is actually beneficial. In situations where data is available on extreme events, distributions with supports outside of the nominal distribution's can lead to over-conservative planning. In particular, this paper deals with application in a disaster response problem in Brazil, where the dataset contains the 2011 Rio de Janeiro floods and mudslides, which is commonly recognized as the worst weather-related natural disaster in Brazilian history. Extending the support beyond a worst-ever disaster situation is not desirable and we therefore omit the Wasserstein distance from consideration in this work.

The main issue with robust stochastic optimization under *f*-divergence is that the practicability of the models vary significantly depending on the choice of the *f*-divergence function. In this paper, we borrow the terminology in Chen et al. (2020), which defines a practicable robust optimization model as one where the computational difficulties of solving the robust counterpart is comparable to that of solving the nominal problem. When first-stage decisions are real-valued in a two-stage linear robust optimization setting, the variation distance yields a linear programming model whereas the Kullback–Leibler function produces a general convex optimization model with exponential terms which, although admitting a self-concordant barrier and being polynomially tractable, is numerically challenging. When first-stage decisions are mixed-integer, the numerical challenge is even more pronounced. For instance, the Kullback–Leibler function produces a mixed-integer model with exponential terms, which is hard to solve to global optimality with state-of-the-art optimization software. Furthermore, many decision-making problems that are posed using event-wise ambiguity require mixed-integer first-stage decisions, such as location-allocation for disaster response.

Conceptually, different *f*-divergences, such as Kullback–Leibler, Burg entropy, and Hellinger distance, are used in robust stochastic optimization because they can capture different attitudes towards ambiguity, or ambiguity aversions. These functions possess properties that decision makers with certain ambiguity aversions would naturally prefer. For instance, the Kullback–Leibler function has a higher second-order derivative than the Hellinger distance, which means that if a decision maker is more averse to large deviations from the nominal distribution, he/she would prefer the Kullback–Leibler over the Hellinger distance. That said, no one can claim that any of these *f*-divergences characterize *exactly* his/her ambiguity aversion, because the latter is not an exactly quantifiable concept.

Given the aforementioned practicability issues in robust stochastic optimization and the understanding that ambiguity aversion is not an exactly quantifiable concept, this paper develops novel divergences that 1) are versatile enough to closely mimic ambiguity aversions conceived by *every* existing *f*-divergence in the literature and 2) maintain a high level of practicability *irrespective* of the existing *f*-divergence they approximate. In addition, our robust stochastic optimization framework remains practicable when first-stage decisions are mixed-integer. In particular, we devise a divergence function based on the infimal convolution of weighted modified variation distances that yields a linear programming robust stochastic optimization counterpart when first-stage decisions are real-valued. Moreover, we derive a piecewise linear divergence that, similar to the infimal convolution, produces a linear programming robust stochastic optimization counterpart, but with greater versatility in modeling ambiguity. Both divergence functions yield mixed-integer programming robust stochastic optimization counterparts when first-stage decisions are mixed-integer. Note that by mimicking *f*-divergences, we mean that our divergences are able to closely approximate the *f*-divergence functions. We make no claim about whether the resulting approximations maintain the information theoretic properties of existing *f*-divergences. The scope of this paper is limited to reproducing the optimal decisions from robust stochastic optimization under *f*-divergences, rather than investigating statistical or information theoretic properties, which, to our knowledge, have not been linked to the decision-making process anywhere in the literature.

However, these functions are limited in that they are non-differentiable. Differentiability is a desirable quality of divergences in some cases, especially when these are used as surrogate loss functions. To create versatile and smooth divergences, we combine the concepts of piecewise linearity and infimal convolution via the Moreau-Yosida regularization. The resulting robust stochastic optimization counterpart is mixed-integer second-order conic at worst.

Our distinct contributions in this paper are summarized as follows:

1. We introduce three divergence functions, based on infimal convolution, piecewise linearity, and Moreau-Yosida regularization. The developed functions offer greater versatility in modeling ambiguity aversions while yielding more practicable robust stochastic optimization counterparts than existing divergence functions. All our resulting formulations can be efficiently solved by commercial optimization packages.
2. We show ways in which our divergence functions can be used to mimic *f*-divergences without losing the superior practicability. This is achieved via least-square fitting for the infimal convolution and the piecewise linear divergences and with a combination of least-square fitting and Monte Carlo integration for the Regularized function. We show that these methods yield closed-form functional approximations of existing *f*-divergence functions.
3. We apply our practicable models to a realistic humanitarian logistics problem that strives to devise an effective and equitable prepositioning relief supply chain based on real Brazilian natural hazards. The academic literature on scenario-based stochastic programming approaches for humanitarian logistics problems has massively increased over the past decades (Sabbaghtorkan et al. 2020). However, most studies still rely on unwarranted assumptions about the scenarios in which the most unrealistic one is perhaps the presumption that the underlying distribution of past disasters is deterministic (Grass and Fischer 2016). This may result in optimal solution prescriptions that are biased towards an inaccurate probability estimation that favors extreme (bad) scenarios and ends up pushing humanitarian assistance in a direction not delivering broader benefits, which is particularly problematic in situations of scarce resources and when people's vulnerability comes into play. Therefore, we propose the first location-allocation humanitarian logistics planning model under ambiguity. Ambiguity has been modelled in other aspects of humanitarian logistics, such as last mile distribution (Zhang et al. 2020) and relief prepositioning (Yang et al. 2021). Our model contains an objective function defined with a utility function and a Gini coefficient, that not only models the effectiveness-equity trade-off in humanitarian operations but also prioritizes allocating scarce resources to more vulnerable areas, based on the well-known FGT poverty measure (Foster et al. 1984). In this case study, we showcase the practicability of our robust stochastic optimization counterparts, the impact of considering ambiguity, and the advantages of our effectiveness-equitable objective function.

The remainder of the paper is organized as follows. Section 2 describes the preliminary concepts and notations of the general robust stochastic optimization problem with *f*-divergences. Section 3 develops the first proposed divergence based on the infimal convolution of modified variation distances. Section 4 develops the piecewise linear divergence. Section 5 presents the Moreau-Yosida regularization of our piecewise linear divergence. Section 6 implements the proposed models in a realistic case study of Brazilian disasters. Section 7 summarizes the work and discusses some future research directions. The appendices to this article contain all formal proofs of theoretical results.

## 2 Robust stochastic optimization framework

Let $\boldsymbol{c} \in \mathbb{R}^{N_1}$ and $\boldsymbol{b} \in \mathbb{R}^{M_1}$ be vectors and $\boldsymbol{A} \in \mathbb{R}^{M_1 \times N_1}$ be a matrix, all of parameters of the first-stage model. The notations $\boldsymbol{d}_\omega \in \mathbb{R}^{N_2}$, $\boldsymbol{D}_\omega \in \mathbb{R}^{M_2 \times N_1}$, $\boldsymbol{E} \in \mathbb{R}^{M_2 \times N_2}$, and $\boldsymbol{f}_\omega \in \mathbb{R}^{M_2}$ represent parameters of the second-stage model, some of which are dependent on scenario $\omega$. Given $S > 1$ and decision variables $\boldsymbol{x}$ and $\boldsymbol{y}$ which are $N_1$-dimensional (consisting of $N_1'$ integers and $N_1''$ continuous, where $N_1' + N_1'' = N_1$) and $N_2$ -dimensional vectors, respectively, and $[S]$ is the set of running indices from 1 to $S$, our robust stochastic optimization framework can be posed as the following problem:

$$
\begin{aligned}
(P1) \quad & \min\ \boldsymbol{c}^\mathsf{T}\boldsymbol{x} + \sup_{p \in \mathcal{P}} \sum_{\omega \in [S]} p_\omega Q(\boldsymbol{x}, \omega) \\
& \text{s.t.} \quad \boldsymbol{A}\boldsymbol{x} \geq \boldsymbol{b} \\
& \qquad\ \ \boldsymbol{x} \in \mathbb{Z}^{N_1'} \times \mathbb{R}^{N_1''} \\
& \text{where} \quad Q(\boldsymbol{x}, \omega) = \min\ \boldsymbol{d}_\omega^\mathsf{T} \boldsymbol{y}_\omega \\
& \text{s.t.} \quad \boldsymbol{D}_\omega \boldsymbol{x} + \boldsymbol{E}\boldsymbol{y}_\omega \geq \boldsymbol{f}_\omega \\
& \qquad\ \ \boldsymbol{y}_\omega \geq \boldsymbol{0},
\end{aligned}
$$

where $\mathcal{P}$ is an ambiguity set used to characterize uncertainties in scenario probabilities. This work positions the framework within the practical context of humanitarian logistics planning. A wide array of preparedness/response planning models follow the robust stochastic optimization paradigm in Model ($P1$), such as location-allocation models, location-prepositioning-allocation models, last-mile distribution models and mass casualty response planning models. The paradigm's main enabling feature is that it allows two-stage scenario-based planning when scenario probabilities are hard to estimate, and this being without restriction on continuity/integrality of first-stage decisions. Humanitarian logistics problems typically involve decisions made before (preparedness) and after (response) uncertainty realizations, which suitwell two-stage modeling. Also, decisions made in the preparedness phase often involve integral components such as facility location, as well as real-valued components such as inventory prepositioning. In addition, these problems rely on scenario-based data, such as disaster occurrences, for which probabilities are unknown or impossible to estimate (this is especially true for low-frequency high-impact events).

We model uncertainties in scenario probabilities through deviations from a given nominal distribution, such that the $f$-divergence of these deviations does not exceed a given threshold $\Xi$. The $f$-divergence between two probability vectors $\boldsymbol{p}$ and $\boldsymbol{q}$ is defined as

$$
I_\phi(\boldsymbol{p}, \boldsymbol{q}) := \sum_{\omega \in [S]} q_\omega \phi\left(\frac{p_\omega}{q_\omega}\right),
$$

where $\phi(z)$ is a convex $f$-divergence function on $z \geq 0$. The limit definitions of $\phi$ when $q_\omega = 0$ are $0\phi(0/0) = 0$ and $0\phi(a/0) = a \lim_{z\to\infty} \phi(z)/z$. The $f$-divergence function has two main properties: 1) It is convex and 2) $\phi(z) > 0$ for every $z \geq 0$, except at $z = 1$ where $\phi(1) = 0$. Our ambiguity set constructed using $f$-divergence is

$$\mathcal{P} = \left\{ \boldsymbol{p} \in \mathbb{R}_+^{|\omega|} : \sum_{\omega\in[S]} p_\omega = 1, \sum_{\omega\in[S]} q_\omega \phi\left(\frac{p_\omega}{q_\omega}\right) \leq \Xi \right\}.$$

In practice, the model requires almost the same inputs from the decision maker as a classical two-stage stochastic programming model. From the uncertainty characterization perspective, the only additional requirements are an $f$-divergence function and a parameter $\Xi$ to represent the decision maker's aversion to deviations from the initial probability estimates. If the decision maker has reasonable certainty about their probability estimates, they would set $\Xi$ to be low (remember that $\phi(1) = 0$ and $\sum_{\omega\in[S]} q_\omega = 1$, which means that $\Xi = 0$ would ensure zero deviations from the initial probability estimates). Else, $\Xi$ would be set to a high value, where the highest value can be obtained from the optimization model $\max\left\{ \sum_{\omega\in[S]} q_\omega \phi\left(\frac{p_\omega}{q_\omega}\right) : \sum_{\omega\in[S]} p_\omega = 1, \boldsymbol{p} \geq 0 \right\}$.

Throughout this paper, all proposed models will satisfy the following two conditions, borrowing terminologies from Hanasusanto and Kuhn (2018):

**Definition 1** (Relatively complete recourse). Model ($P1$) has relatively complete recourse when, for any given $\omega \in [S]$, $Q(\boldsymbol{x}, \omega)$ is feasible for all $\boldsymbol{x} \in \{\mathcal{X} : \boldsymbol{A}\boldsymbol{x} \geq \boldsymbol{b}\}$.

**Definition 2** (Sufficiently expensive recourse). Model ($P1$) has sufficiently expensive recourse when, for any given $\omega \in [S]$, the dual of $Q(\boldsymbol{x}, \omega)$ is feasible.

If Model ($P1$) has both relatively complete recourse and sufficiently expensive recourse, $Q(\boldsymbol{x}, \omega)$ is feasible and finite. These two conditions are not overly restrictive and are satisfied by many problems or can be enforced by induced constraints. For instance, in humanitarian problems where demands have to be satisfied, inclusion of unmet demands with finite cost penalties can ensure feasibility and finiteness of $Q(\boldsymbol{x}, \omega)$. Another possibility is to include first-stage constraints that ensure enough facilities (propositioned inventory) are located (deployed) to satisfy all demands. The following theorem, adapted from Bayraksan and Love (2015), shows that Model ($P1$) can be reformulated in terms of the so-called convex conjugate of $\phi$. The convex conjugate, denoted by the function $\phi^* : \mathbb{R} \to \mathbb{R} \cup \{\infty\}$, is defined as

$$\phi^*(s) = \sup_{z\geq 0}\{sz - \phi(z)\}.$$

**Theorem 1** (Bayraksan and Love 2015) *Model* ($P1$) *is equivalent to*

$$\begin{aligned} \min \quad & \boldsymbol{c}^\top \boldsymbol{x} + \lambda \Xi + \mu^+ - \mu^- + \lambda \sum_{\omega \in [S]} q_\omega \phi^* \left( \frac{\boldsymbol{d}_\omega^\top \boldsymbol{y}_\omega - \mu^+ + \mu^-}{\lambda} \right) \\ \text{s.t.} \quad & \boldsymbol{A}\boldsymbol{x} \ge \boldsymbol{b} \\ & \boldsymbol{D}_\omega \boldsymbol{x} + \boldsymbol{E}\boldsymbol{y}_\omega \ge \boldsymbol{f}_\omega \quad \forall \omega \in [S] \\ & \boldsymbol{x} \in \mathbb{Z}^{N_1'} \times \mathbb{R}^{N_1''} \\ & \boldsymbol{y}_\omega \ge \boldsymbol{0} \quad \forall \omega \in [S] \\ & \lambda, \mu^+, \mu^- \ge 0. \end{aligned}$$

Depending on the *f*-divergence used, the robust counterpart in Theorem 1 will have different levels of practicability. For example, with the Hellinger distance, the model is mixed-integer conic quadratic, whereas the Kullback–Leibler divergence leads to a general mixed-integer optimization model containing exponential terms. The only *f*-divergence known in the literature that yields a mixed-integer robust counterpart for Model ($P1$) is the variation distance. However, the variation distance is linear and lack versatility in portraying ambiguity preferences. For reference, we provide Table 1 that lists some popular examples of *f*-divergences for cases where closed-form conjugates exist. Humanitarian logistics problems are varied, with contrasting levels of past data consistency. For instance, disaster impacts on a very disaster-prone area with stable infrastructure would likely resemble impacts of similar past disasters, which would make probabilities reliably estimable from past data. In other cases, either events have rare features (e.g. COVID-19 pandemic), or infrastructures change significantly, meaning that disaster impacts vary unpredictably, thereby making probability estimates unreliable or unrealistic. The contrasting

**Table 1** Examples of *f*-divergences with their conjugates and the tractability of their robust stochastic optimization counterparts

| Divergence | $\phi(z), z \ge 0$ | $I_\phi(\boldsymbol{p}, \boldsymbol{q})$ | $\phi^*(s)$ | Tractability of robust stochastic optimization counterpart |
|---|---|---|---|---|
| Kullback–Leibler | $z \log z - z + 1$ | $\sum p_\omega \log(\frac{p_\omega}{q_\omega})$ | $e^s - 1$ | Mixed-integer with exponential terms |
| Burg entropy | $-\log z + z - 1$ | $\sum q_\omega \log(\frac{q_\omega}{p_\omega})$ | $-\log(1-s), s < 1$ | Mixed-integer with logarithmic terms |
| $\chi^2$-distance | $\frac{1}{z}(z-1)^2$ | $\sum \frac{(p_\omega - q_\omega)^2}{p_\omega}$ | $2 - 2\sqrt{1-s}, s < 1$ | Mixed-integer conic quadratic |
| Hellinger distance | $(\sqrt{z} - 1)^2$ | $\sum(\sqrt{p_\omega} - \sqrt{q_\omega})^2$ | $\frac{s}{1-s}, s < 1$ | Mixed-integer conic quadratic |
| Variation distance | $\lvert z - 1 \rvert$ | $\sum \lvert p_\omega - q_\omega \rvert$ | $\begin{cases} -1 & s \le -1 \\ s & -1 \le s \le 1 \end{cases}$ | Mixed-integer |

Adapted from (Ben-Tal et al. 2013)

levels of past data consistency make it necessary to offer a wide variety of functions to model different ambiguity preferences and thus capture a wide range of realities via robust stochastic modeling. Moreover, these models tend to be large-scale, making it necessary to find ways to solve them within reasonable times.

## 3 A divergence construction using infimal convolution

While possessing desirable properties such as convexity and versatility in modeling diverse ambiguity aversions, *f*-divergences are never *exact* measures of ambiguity aversion because the latter is an inherently subjective concept. Even in cases of consistent disaster impacts on stable infrastructure, there still exists a personal element in ambiguity aversion that is attached to each decision maker's preferences (political ideologies of different governments, *modus operandi* of different humanitarian organizations) and is not dependent on the impacted area *per se*. In this section and subsequent ones, we seek to propose new divergence functions that offer greater versatility than *f*-divergences, mainly because 1) they can be built based on the ambiguity preferences of the decision maker and 2) they can mimic existing *f*-divergences while being more practicable.

An infimal convolution of $\mathcal{D}$ closed convex functions $\psi_1, \dots, \psi_{\mathcal{D}}$ is defined by

$$(\psi_1 \square \cdots \square \psi_{\mathcal{D}})(z) := \inf_{\sum_{d\in[\mathcal{D}]} r_d = z} \left\{ \sum_{d\in[\mathcal{D}]} \psi_d(r_d) \right\}.$$

The following theorem shows that under specific modifications, the infimal convolution of the $\mathcal{D}$ modified variation distances is itself a divergence function, which we define to be a function that is convex and positive everywhere on the domain, except at $z = 1$, where it is equal to zero. This definition follows that of general statistical divergences, with the addition of convexity, which is essential in optimization modeling. The theorem also shows that under further conditions, the infimal convolution is equivalent to the variation distance.

**Theorem 2** *For any* $\omega \in [S]$, *the infimal convolution* $\inf_{\sum_{d\in[\mathcal{D}]} r_d = p_\omega/q_\omega} \{\sum_{d\in[\mathcal{D}]} \pi_d \mid \mathcal{D} r_d - 1 \mid\}$, *where* $\boldsymbol{\pi} > \mathbf{0}$, *is* 1) *convex*, 2) *equal to zero when* $p_\omega = q_\omega$ *and* 3) *positive when* $p_\omega \neq q_\omega$. *Furthermore, when* $\pi_d = \frac{1}{\mathcal{D}}$, $\forall d \in [\mathcal{D}]$, *this infimal convolution is equivalent to the variation distance* $\mid \frac{p_\omega}{q_\omega} - 1 \mid$.

Therefore, the infimal convolution $\inf_{\sum_{d\in[\mathcal{D}]} r_d = p_\omega/q_\omega} \{\sum_{d\in[\mathcal{D}]} \pi_d \mid \mathcal{D} r_d - 1 \mid\}$ is a divergence measure, irrespective of the choice of positive $\boldsymbol{\pi}$. The decision maker is therefore free to choose values of $\boldsymbol{\pi}$ such that the resulting infimal convolution matches his/her ambiguity aversion. Independent of the choice, the solvable form of Model ($P1$) is a mixed-integer programming model, as formally shown in the next theorem.

**Theorem 3** *Under the divergence* $\inf_{\sum_{d\in[\mathcal{D}]} r_d = p_\omega/q_\omega} \{\sum_{d\in[\mathcal{D}]} \pi_d \mid \mathcal{D} r_d - 1 \mid\}$, $\forall \omega \in [S]$, *Model* ($P1$) *is equivalent to the following mixed-integer programming model*:

$$
\begin{aligned}
(P2)\ \min \quad & \boldsymbol{c}^{\mathsf{T}}\boldsymbol{x} + \lambda\Xi + \mu^{+} - \mu^{-} + \sum_{\substack{\omega \in [S] \\ d \in [\mathcal{D}]}} q_{\omega} z_{\omega d} \\
\text{s.t.} \quad & \boldsymbol{A}\boldsymbol{x} \geq \boldsymbol{b} \\
& \boldsymbol{D}_{\omega}\boldsymbol{x} + \boldsymbol{E}\boldsymbol{y}_{\omega} \geq \boldsymbol{f}_{\omega} \quad \forall \omega \in [S] \\
& z_{\omega d} \geq -\lambda \pi_{d} \quad \forall \omega \in [S], d \in [\mathcal{D}] \\
& D z_{\omega d} \geq \boldsymbol{d}_{\omega}^{\mathsf{T}}\boldsymbol{y}_{\omega} - \mu^{+} + \mu^{-} \quad \forall \omega \in [S], d \in [\mathcal{D}] \\
& \boldsymbol{d}_{\omega}^{\mathsf{T}}\boldsymbol{y}_{\omega} - \mu^{+} + \mu^{-} \leq \lambda D \pi_{d} \quad \forall \omega \in [S], d \in [\mathcal{D}] \\
& \boldsymbol{x} \in \mathbb{Z}^{N_1'} \times \mathbb{R}^{N_1''} \\
& \boldsymbol{y}_{\omega} \geq \boldsymbol{0}, z_{\omega d} \in \mathbb{R} \quad \forall \omega \in [S], d \in [\mathcal{D}] \\
& \lambda, \mu^{+}, \mu^{-} \geq 0 \quad \forall \omega \in [S].
\end{aligned}
$$

Whilst different values of $\boldsymbol{\pi}$ yield different divergence functions, the practicability of the robust stochastic optimization counterpart is unchanged. The vector $\boldsymbol{\pi}$ can thus be used to tailor the divergence function to match the decision maker's ambiguity aversion. In the absence of a clear preference, $\boldsymbol{\pi}$ can be computed by using an existing $f$-divergence function as a reference. We hereafter provide decision makers with a method to estimate $\boldsymbol{\pi}$ based on how well it makes the infimal convolution of modified variation distances (ICV) fit an existing $f$-divergence function of choice. This helps reduce the practical arbitrariness involved in choosing $\boldsymbol{\pi}$ and thus, in finding an appropriate divergence function. In addition, our method does not affect the practicability of model implementation, since it allows $\boldsymbol{\pi}$ to be computed without solving any auxiliary model. With the resulting value, decision makers can easily plot and visually assess the resulting divergence function. The next theorem shows that if the ICV is a least-square fit of an existing $f$-divergence function, $\boldsymbol{\pi}$ has a closed-form optimal value. To prove that, we first define the sum of square of differences (SSD) for the ICV as

$$
\mathrm{SSD}(\boldsymbol{\pi}) := \int_0^H \left( \inf_{\sum_{d \in [\mathcal{D}]} r_d = z} \left\{ \sum_{d \in [\mathcal{D}]} \pi_d \mid \mathcal{D} r_d - 1 \mid \right\} - \phi(z) \right)^2 \mathrm{d}z,
$$

where $H = \max_{\omega \in [S]} \{ \frac{\bar{q}_{\omega}}{q_{\omega}} \}$ and $\bar{q}_{\omega}$ is an upper bound on the probability of scenario $\omega$.

**Theorem 4** *The minimizer* $\boldsymbol{\pi}^*$ *of* SSD *is such that*

$$
\pi^*_{\underline{d}} = \frac{3\Phi}{\mathcal{D}((H-1)^3 + 1)},
$$

*where* $\Phi = \int_0^H \phi(z) \mid z - 1 \mid \mathrm{d}z$, *for any one arbitrary* $\underline{d} \in [\mathcal{D}]$ *and* $\pi_d^*$ *is any value greater than* $\pi^*_{\underline{d}}$ *for all* $d \in [\mathcal{D}] \backslash \{\underline{d}\}$.

Interestingly, a further practicability improvement can be achieved with the least-square ICV (LS-ICV), i.e., the ICV such that the value of $\boldsymbol{\pi}$ minimizes the SSD. This is because of the following lemma, which proves that the LS-ICV of any $f$-divergence function $\phi$ is simply a weighted variation distance that is independent of $\mathcal{D}$. This follows from the fact that under the least-square weights, the minimum SSD is independent of $\mathcal{D}$.

**Lemma 5** *The LS-ICV of* $\phi(z)$ *is* $\dfrac{3\Phi}{(H-1)^3+1} \mid z-1 \mid$.

Because of Lemma 5, the solvable form of Model ($P1$) under the LS-ICV is independent of the value of $\mathcal{D}$, which means that it admits fewer variables and constraints than Model ($P2$), as shown by the following theorem.

**Theorem 6** *Under the divergence* $\dfrac{3\Phi}{(H-1)^3+1} \mid \dfrac{p_\omega}{q_\omega} - 1 \mid, \forall \omega \in [S]$, *Model* ($P1$) *is equivalent to the following mixed-integer programming model*:

$$
\begin{aligned}
(P3)\ \min\ & \boldsymbol{c}^{\mathsf{T}}\boldsymbol{x} + \lambda\Xi + \mu^+ - \mu^- + \sum_{\omega\in[S]} q_\omega z_\omega \\
\text{s.t.}\ & \boldsymbol{A}\boldsymbol{x} \geq \boldsymbol{b} \\
& \boldsymbol{D}_\omega \boldsymbol{x} + \boldsymbol{E}\boldsymbol{y}_\omega \geq \boldsymbol{f}_\omega \quad \forall \omega \in [S] \\
& z_\omega \geq -\frac{3\lambda\Phi}{(H-1)^3+1} \quad \forall \omega \in [S] \\
& z_\omega \geq \boldsymbol{d}_\omega^{\mathsf{T}}\boldsymbol{y}_\omega - \mu^+ + \mu^- \quad \forall \omega \in [S] \\
& \boldsymbol{d}_\omega^{\mathsf{T}}\boldsymbol{y}_\omega - \mu^+ + \mu^- \leq \frac{3\lambda\Phi}{(H-1)^3+1} \quad \forall \omega \in [S] \\
& \boldsymbol{x} \in \mathbb{Z}^{N_1'} \times \mathbb{R}^{N_1''} \\
& \boldsymbol{y}_\omega \geq \boldsymbol{0}, z_\omega \in \mathbb{R} \quad \forall \omega \in [S] \\
& \lambda, \mu^+, \mu^- \geq 0 \quad \forall \omega \in [S].
\end{aligned}
$$

Model ($P1$) is therefore shown to be a mixed-integer programming model under the LS-ICV, with fewer variables and constraints than under the general ICV. This means that when an existing $f$-divergence function is used as a reference to estimate optimal weights in the infimal convolution using the least-square method, the practicability improves. In a way, by offering decision makers a method to tailor the new divergence function according to existing ones, we have also offered them a way to improve the solvability of the resulting decision-making problem. Equipped with this interesting set of results, in subsequent sections, we provide decision makers with new divergence functions that are even more versatile than the LS-ICV and that can capture the behaviors of existing $f$-divergence functions even more closely, while keeping decision-making models easily solvable and least-square estimates of parameters computable in closed forms.

## 4 A piecewise linear divergence

To offer further improve versatility in ambiguity aversion modeling, with little sacrifice to practicability, we propose a piecewise linear divergence. This divergence allows the decision maker to portray different ambiguity aversions in different ranges of probability deviations. The piecewise linear divergence with $P$ pieces is defined by $\max_{p\in[P]}\{a_p z - b_p\}$, such that $\phi(1) = 0$ and $\phi(z) > 0$ for all other $z$. In the same vein as in the infimal convolution case, the decision maker can choose values of $\boldsymbol{a}$ and $\boldsymbol{b}$ that better capture his/her ambiguity aversion. The difference here is that the choices can be varied according to the probability deviation range, which means that the decision maker is able to tailor his/her attitude towards ambiguity depending on ranges of probabilities and on how far from nominal the probability is. The piecewise linear divergence has the added advantage of being able to mimic more closely the behaviours of existing $f$-divergences, while offering better practicability.

Piecewise linear fitting is a challenging task that is generally solved heuristically. There exists a whole body of literature on piecewise linear fitting and the main methods proposed rely on mixed-integer programming approaches (Rebennack and Krasko 2020; Toriello and Vielma 2012) or non-convex real-valued approaches (Rebennack and Kallrath 2015). In our case, such methods would involve solving additional models (non-convex models may even require additional algorithmic developments), which would reduce the practicability of our robust stochastic optimization framework. Indeed, this work focuses on producing robust counterparts with the same level of numerical difficulty as the nominal problem,i.e., practicable robust counterparts. To ensure easily computable piecewise linear approximations, we impose restrictions on the positions of breakpoints. Moreover, we show in the following proposition that a piecewise linear divergence that is a piecewise least-square fit of an existing $f$-divergence, which we term LS-PL, can be expressed as a closed-form recursive function that does not require the solution of additional models. We also show in Theorem 10 that our function can be dualized to produce a mixed-integer programming robust counterpart.

**Proposition 7** *Suppose a piecewise linear divergence is constructed such that fitting is done separately for ranges $z \le 1$ and $1 \le z \le H$. If $L$ and $U$ pieces with equally spaced intersection points are used for $z \le 1$ and $1 \le z \le H$, respectively, the piecewise linear divergence that fits an existing divergence $\phi(z)$, such that pieces are sequentially least-square fitted starting at $z = 1$, is given by*

$$\mathcal{G}(z) = \begin{cases} \max_{l\in[L]} \left\{ \left(3L^3\Psi_l^a - \frac{3L}{2} f_{l+1}^a\left(\frac{l}{L}\right)\right)\left(\frac{l}{L} - z\right) + f_{l+1}^a\left(\frac{l}{L}\right)\right\} & \text{if } 0 \le z \le 1, \\ \max_{u\in[U]} \left\{ \left(\frac{3}{\Delta^3}\Psi_u^b - \frac{3}{2\Delta} f_{u-1}^b(1+(u-1)\Delta)\right)(z - 1 - (u-1)\Delta) \right. & \\ \left. + f_{u-1}^b(1+(u-1)\Delta)\right\} \quad \text{if } 1 \le z \le H, & \end{cases},$$

*where* $f^a_{L+1}\Big(\frac{l}{L}\Big) = 0$, $\Psi^a_l = \int_{(l-1)/L}^{l/L} \phi(z)\Big(\frac{l}{L} - z\Big)\mathrm{d}z$, $\Delta = \frac{H-1}{U}$, $f^b_0(1+(u-1)\Delta) = 0$ *and* $\Psi^b_u = \int_{1+(u-1)\Delta}^{1+u\Delta} \phi(z)(z-1-(u-1)\Delta)\mathrm{d}z$.

This shows that under mild conditions, to obtain the least-square piecewise linear fit of an existing *f*-divergence function, we need not solve the complicated definite integral in the SSD formula. In addition, the SSD of the fit can be calculated in closed form, as indicated in the following corollary.

**Corollary 8** *The minimum SSD when the piecewise linear function $\mathcal{G}(z)$ is used to fit an existing divergence $\phi(z)$ is*

$$
\begin{aligned}
\mathrm{SSD}^* = {} & \sum_{p\in[L]} \frac{\Big(a^*_p \frac{p}{L} - b^*_p - \phi(\frac{p}{L})\Big)^3}{3\Big(a^*_p - \phi'(\frac{p}{L})\Big)} - \sum_{p\in[L]} \frac{\Big(a^*_p \frac{p-1}{L} - b^*_p - \phi(\frac{p-1}{L})\Big)^3}{3\Big(a^*_p - \phi'(\frac{p-1}{L})\Big)} \\
& + \sum_{p\in[L+U]\setminus[L]} \frac{\Big(a^*_p(1+(p-L)\Delta) - b^*_p - \phi(1+(p-L)\Delta)\Big)^3}{3\Big(a^*_p - \phi'(1+(p-L)\Delta)\Big)} \\
& - \sum_{p\in[L+U]\setminus[L]} \frac{\Big(a^*_p(1+(p-L-1)\Delta) - b^*_p - \phi(1+(p-L-1)\Delta)\Big)^3}{3\Big(a^*_p - \phi'(1+(p-L-1)\Delta)\Big)},
\end{aligned}
$$

*where*

$$
a^*_p = \begin{cases} -\Big(3L^3\Psi^a_p - \frac{3L}{2} f^a_{p+1}\Big(\frac{p}{L}\Big)\Big) & \text{if } p \in [L] \\ \Big(\frac{3}{\Delta^3}\Psi^b_{p-L} - \frac{3}{2\Delta} f^b_{p-L-1}(1+(p-L-1)\Delta\Big) & \text{if } p \in [L+U]\setminus[L], \end{cases}
$$

*and*

$$
b^*_p = \begin{cases} \Big(\frac{3p}{2} - 1\Big) f^a_{p+1}\Big(\frac{p}{L}\Big) - 3pL^2\Psi^a_p & \text{if } p \in [L] \\ \frac{3+(p-L-1)\Delta}{\Delta^3}\Psi^b_{p-L} - \Big(1 + \frac{3}{2\Delta} + \frac{3(p-L-1)}{2}\Big) f^b_{p-L-1}(1+(p-L-1)\Delta) \\ \text{if } p \in [L+U]\setminus[L] \end{cases}
$$

*and* $f^a_{L+1}\Big(\frac{p}{L}\Big) = 0$, $\Psi^a_p = \int_{(p-1)/L}^{p/L} \phi(z)\Big(\frac{p}{L} - z\Big)\mathrm{d}z$, $\Delta = \frac{H-1}{U}$, $f^b_0(1+(p-L-1)\Delta) = 0$, $\Psi^b_{p-L} = \int_{1+(p-L-1)\Delta}^{1+(p-L)\Delta} \phi(z)(z-1-(p-L-1)\Delta)\mathrm{d}z$ *and* $\phi'(m)$ *is the derivative of* $\phi$ *evaluated at point m.*

The following theorem states that the LS-PL is always at least as good as the LS-ICV in mimicking the behaviours of existing *f*-divergences.

**Theorem 9** *The SSD of the LS-PL is always less than or equal to that of the LS-ICV.*

Finally, we show that, similar to the LS-ICV, the LS-PL also yields a mixed-integer programming model

**Theorem 10** *Under the divergence* $\mathcal{G}\left(\frac{p_\omega}{q_\omega}\right)$, *Model* $(P1)$ *is equivalent to the following mixed-integer programming model*:

$$
\begin{aligned}
(P4)\ \min\quad & \boldsymbol{c}^\top \boldsymbol{x} + \lambda \Xi + \mu^+ - \mu^- + \sum_{\omega \in [S]} q_\omega z_\omega \\
\text{s.t.}\quad & \boldsymbol{Ax} \geq \boldsymbol{b} \\
& \boldsymbol{D}_\omega \boldsymbol{x} + \boldsymbol{E}\boldsymbol{y}_\omega \geq \boldsymbol{f}_\omega \quad \forall \omega \in [S] \\
& z_\omega \geq \frac{b_p^* - b_{p+1}^*}{(a_p^* - a_{p+1}^*)}(\boldsymbol{d}_\omega^\top \boldsymbol{y}_\omega - \mu^+ + \mu^- - \lambda a_p^*) + \lambda b_p^* \quad \forall \omega \in [S], p \in [L+U-1] \\
& \boldsymbol{x} \in \mathbb{Z}^{N_1'} \times \mathbb{R}^{N_1''} \\
& \boldsymbol{y}_\omega \geq \boldsymbol{0}, z_\omega \in \mathbb{R} \quad \forall \omega \in [S] \\
& \lambda, \mu^+, \mu^- \geq 0.
\end{aligned}
$$

Therefore, the piecewise linear divergence is at least as good as the infimal convolution divergence at replicating ambiguity aversions modeled through existing *f*-divergences, while yielding a solvable form for Model $(P2)$ that is close in practicability to the solvable form under the infimal convolution divergence. Notice that $\boldsymbol{a}^*$ and $\boldsymbol{b}^*$ are calculated a priori and could have easily been replaced with pre-defined values if the decision maker wishes to model his/her ambiguity aversion differently from existing *f*-divergences.

## 5 Convolution with Moreau-Yosida regularization

One issue in the divergences proposed so far is the lack of differentiability. Our infimal convolution divergence is not differentiable at the point where it is equal to zero and the piecewise linear divergence is not differentiable at the points of intersection of pieces. Differentiability of *f*-divergences is especially required in the machine learning context (Bartlett et al. 2006) where these divergences are used as surrogate loss functions. In the humanitarian logistics context, differentiability allows decision makers to understand ambiguity aversion based on rates of change of divergence functions with probability deviations. In a way, the divergence function value is a penalty for deviations from nominal probabilities. For instance, if a disaster is consistent with historical data and the probability can be estimated with high level of confidence, decision makers would choose divergence functions with high rates of change, so as to prevent large probability deviations. In the piecewise linear case, decision makers know the rates of change in intervals between intersecting pieces, but not at points of intersection. This becomes problematic if the number of pieces is large, as non-differentiability spreads across a large number of points. The class of $f$-divergence functions shown in (Ben-Tal et al. 2013) contains a mix of differentiable and non-differentiable functions. In the spirit of proposing divergences that have widespread applicability, we combine the infimal convolution with piecewise linearity to produce smooth and differentiable piecewise linear divergences that can mimic the behaviours of existing $f$-divergences. The smooth piecewise linear divergence can more accurately mimic existing smooth $f$-divergences. We use the Moreau-Yosida regularization, itself defined as an infimal convolution, due to its desirable conjugacy and closed-form property.

**Definition 3** The Moreau-Yosida regularization $\mathcal{Y}$ of a closed convex function $g$ is

$$\mathcal{Y}(z) := \min_{s \in \mathbb{R}^n} \left\{ g(s) + \frac{1}{2} \|z - s\|_M^2 \right\},$$

where $\|z - s\|_M^2 = (z - s)^T M(z - s)$, $a \in \mathbb{R}^n$, and $M$ is a symmetric positive definite $n \times n$ matrix.

The Moreau-Yosida regularization of $\mathcal{G}$ (which is a closed convex function) is $\mathcal{Y}(z) = \min_{s \in \mathbb{R}} \{\max_{p \in [L+U]} \{a_p^* s - b_p^*\} + \frac{m}{2}(z - s)^2\}$, where $m$ is a positive scalar. This regularization smoothes the piecewise linear divergence and makes it differentiable. For a comprehensive exposition of the properties of the Moreau-Yosida regularization, we refer the reader to (Lemaréchal and Sagastizábal 1997).

The regularization of $\mathcal{G}$, explicitly stated, is:

$$\mathcal{Y}(z) = \begin{cases} a_1^* z - b_1^* - \dfrac{(a_1^*)^2}{2m} & \mathit{if}\ l_0 + \dfrac{a_1^*}{m} \le z < l_1 + \dfrac{a_1^*}{m} \\ a_p^* z - b_p^* - \dfrac{(a_p^*)^2}{2m} & \mathit{if}\ l_{p-1} + \dfrac{a_p^*}{m} < z < l_p + \dfrac{a_p^*}{m}, p \in \{2, \dots, L+U-1\} \\ a_{L+U}^* z - b_{L+U}^* - \dfrac{(a_{L+U}^*)^2}{2m} & \mathit{if}\ l_{L+U-1} + \dfrac{a_{L+U}^*}{m} < z \le l_{L+U} + \dfrac{a_{L+U}^*}{m} \\ \min_{p \in [L+U-1]} \{a_p^* l_p - b_p^* - \frac{m}{2}(z - l_p)^2\} & \mathit{otherwise}, \end{cases}$$

where $l_0 = -\frac{a_1^*}{m}$, $l_{L+U} = H - \frac{a_{L+U}^*}{m}$ and the remaining $l_p$ are the intersection points of pieces in the LS-PL. The above formulation follows from differentiating $\mathcal{Y}$, setting it to zero in regions of the domain where $\mathcal{G}$ is differentiable and taking the formulation at intersection points where it is not. The value of $m$ affects the accuracy with which the regularized fit can mimic the behaviours of existing $f$-divergence functions. In the following proposition, we derive the value of $m$ that makes $\mathcal{Y}$ a least-square fit of $\phi$ using Monte Carlo integration over a specific subset of the domain.

**Proposition 11** *If $\mathcal{Z}$ points are sampled uniformly from the range $[l'_{p-1}+\epsilon, l'_p-\epsilon]$, $\forall p \in [L+U]$, where $l'_0 = 0$, $l'_{L+U} = H$ and $l'_p = l_p$ for $p \in \{2,\dots,L+U-1\}$, the value of m that makes $\mathcal{Y}$ a least-square fit of $\phi$ is obtained from the following quadratic programming problem*

$$m^* = \arg\min_{m>0}\left\{\sum_{p\in[L+U]}\frac{l'_p-l'_{p-1}-2\epsilon}{\mathcal{Z}}\sum_{i\in[\mathcal{Z}]}\left(a^*_p z_{ip}-b^*_p-\frac{(a^*_p)^2}{2m}-\phi(z_{ip})\right)^2\right\},$$

*where* $\epsilon \geq \dfrac{\max_{p\in[L+U]}\{a^*_p\}}{m^*}$.

Without running the optimization model in the proposition, the value of $m^*$ can be obtained in closed form by setting the derivative to zero for cases where the critical point of SSD is a minimizer. For such cases,

$$m^* = \frac{\mathcal{Z}\sum_{p\in[L+U]}(l'_p-l'_{p-1}-2\epsilon)(a^*_p)^4}{2\left(\sum_{p\in[L+U]}(l'_p-l'_{p-1}-2\epsilon)(a^*_p)^2\sum_{i\in[\mathcal{Z}]}(a^*_p z_{ip}-\phi(z_{ip})-b^*_p)\right)}.$$

The value of $\epsilon$ seems to be determined a posteriori in this proposition, which is problematic since it is needed in the calculation of $m^*$. A simple solution to this is to initialize $\epsilon$ and then iteratively lower it until before $\epsilon \geq \dfrac{\max_{p\in[L+U]}\{a^*_p\}}{m^*}$ is violated. With the regularized divergence, Model $(P1)$ is solvable as a mixed-integer second order conic problem, as shown in the following theorem.

**Theorem 12** *Under the divergence $\mathcal{Y}(\frac{p_\omega}{q_\omega})$, Model $(P1)$ is equivalent to the following mixed-integer second order conic problem*:

$$\begin{aligned}
(P6)\ \min\quad & \boldsymbol{c}^\top\boldsymbol{x}+\lambda\Xi+\mu^+-\mu^-+\sum_{\omega\in[S]}q_\omega\left(z_\omega+\frac{1}{2m}\tau_\omega\right)\\
\text{s.t.}\quad & \boldsymbol{A}\boldsymbol{x}\geq\boldsymbol{b}\\
& \boldsymbol{D}_\omega\boldsymbol{x}+\boldsymbol{E}\boldsymbol{y}_\omega\geq\boldsymbol{f}_\omega \quad \forall\omega\in[S]\\
& z_\omega \geq \frac{b^*_p-b^*_{p+1}}{(a^*_p-a^*_{p+1})}(\boldsymbol{d}_\omega^\top\boldsymbol{y}_\omega-\mu^++\mu^--\lambda a^*_p)+\lambda b^*_p \quad \forall\omega\in[S], p\in[L+U-1]\\
& \alpha_\omega \geq \boldsymbol{d}_\omega^\top\boldsymbol{y}_\omega-\mu^++\mu^- \quad \forall\omega\in[S]\\
& \tau_\omega+\lambda \geq \left\|\begin{pmatrix}\sqrt{2}\alpha_\omega\\ \tau_\omega\\ \lambda\end{pmatrix}\right\|_2 \quad \forall\omega\in[S]\\
& \boldsymbol{x}\in\mathbb{Z}^{N'_1}\times\mathbb{R}^{N''_1}\\
& \boldsymbol{y}_\omega\geq\boldsymbol{0}, \alpha_\omega, \tau_\omega\geq 0, z_\omega\in\mathbb{R} \quad \forall\omega\in[S]\\
& \lambda, \mu^+, \mu^-\geq 0.
\end{aligned}$$

The price of differentiability, when it is administered via the Moreau-Yosida regularization is, then, a loss of tractability. A point that is worth noting is that for cases such as the Kullback–Leibler divergence and the Burg Entropy, the regularization can replicate closely the ambiguity aversions they represent while actually improving the practicability from a general mixed-integer problem with exponential/logarithmic terms to a mixed-integer second order conic problem, which can be very desirable if running time, approximation closeness, and differentiability are crucial in implementing the models in practical decision support systems.

A crucial point to note about the LS-ICV, LS-PL and the regularized LS-PL is that they can also mimic $f$-divergences that do not have closed-form convex conjugates, such as the $J$-divergence (see (Ben-Tal et al. 2013)). Within an even bigger picture, our proposed divergences allow the decision maker to altogether avoid the tedious derivation of convex conjugates, and thus use a greater variety of divergence measures for which the conjugates are not readily available, such as Bregman divergences and $\alpha$-divergences.

## 6 Numerical study: the Brazilian humanitarian supply chain

In July 2013, the Brazilian National Department for Civil Protection and Defense (SEDEC) reached an agreement with the Brazilian Postal Office Service (‘*Correios*’) to preposition relief supplies, such as food, water, and medicine, among others, in municipalities in major states in Brazil. In case of a disaster, SEDEC would transship these supplies between municipalities depending on victim needs, with the local civil defenses thereafter taking over last-mile distributions. In an unexpected turn of events, the endeavour was abandoned due to the lack of coordination of logistics operations and the high implementation costs. SEDEC instead resorted to a supplier selection strategy based on the ability to fulfil relief demands within 192 h for the North region and 96 h for the remaining regions. This strategy is flawed since most relief supplies are needed within the first 48 critical hours after a disaster has struck. In addition, this strategy is blind to imbalances in socioeconomic vulnerabilities, which is a prominent reality in Brazil. Many recent disasters in Brazil are actually a consequence of social processes whose main characteristic is the unequal distribution of opportunities or social inequality that push poorer people to risky areas or to informal settlements often placed in slopes and floodplains that lack basic infrastructure (Carmo and Anazawa 2014; Alem et al. 2021).

Considering that prepositioning of relief aid is one of the most effective preparedness strategies to deal with most types of disasters, but acknowledging it can be prohibitively complicated and expensive, we propose to build a location-allocation network for strategic relief aid prepositioning focused on the protection of the most vulnerable groups in a disaster aftermath. Our approach is meant to be implemented by an entity such as the Federal Government, who will be in charge of allocating relief aid to help the affected municipalities within the first 48 critical hours. When the lack of resources makes it impossible to maintain high levels of prepositioned stock to supply all incoming aid requests at once, our formulation will judiciously select locations that should be prioritized according to their *utility* profile, which includes a measure of vulnerability, amongst others.

Formally stated, our problem concerns a given geographical region prone to natural hazards and is defined with respect to a number of settlements (towns, villages or even neighbourhoods), each of which we refer to as an affected area, that exhibit different characteristics such as demographic and socioeconomic profiles. In a disaster aftermath, these affected areas may require relief aid from prepositioned stockpiles of supplies in response facilities. Humanitarian logisticians must decide here-and-now on where to set up response facilities and to what level critical supplies must be stored in these facilities. After a disaster strikes, they must make wait-and-see decisions on how to assign victim needs to the set-up response facilities.

The proposed optimization model uses the following notation. $\mathcal{N}$ is the set of potential response facility locations, $\mathcal{A}$ is the set of affected areas, $\mathcal{L}$ is the set of response facility sizes, $\mathcal{R}$ is the set of relief aid supplies, and $[S]$ is the set of indices running from 1 to $S$, in which $S$ denotes the number of disaster scenarios. We denote by $c^{\mathrm{o}}_{\ell n}$ the fixed cost of setting up response facility of size $\ell$ at location $n$, $c^{\mathrm{p}}_{rn}$ the prepositioning cost of relief aid $r$ at location $n$, and $c^{\mathrm{d}}_{an}$ the unit cost of shipping relief aid supplies from location $n$ to affected area $a$. The capacity (in volume) of response facility of size $\ell$ is given by $\kappa^{\mathrm{resp}}_{\ell}$, and the units of storage space required by relief aid $r$ is given by $\rho_r$. Prepositioning of relief aid $r$ at location $n$ implies a minimum quantity of $\theta^{\min}_{rn}$, and the overall quantity of relief aid $r$ to be stockpiled is denoted by $\theta^{\max}_{r}$. Pre- and post-disaster humanitarian operations are budgeted at $\eta$ and $\eta'$, respectively. Victim needs for relief aid $r$ at affected area $a$ in scenario $\omega$ are represented by $d_{ra\omega}$ and the nominal probability of occurrence of scenario $\omega$ is given by $p_{\omega}$. Finally, we denote by $u_{ran\omega}$ the utility of assignment of 100% of relief aid $r$ of affected area $a$ to response facility $n$ in scenario $\omega$. Our model entails the following decision variables: $P_{rn}$ is the quantity of relief aid $r$ prepositioned at response facility $n$, $Y_{\ell n}$ is the binary variable that indicates whether or not the response facility of category $\ell$ is established at location $n$, $X_{ran\omega}$ is the fraction of relief aid $r$ of affected area $a$ satisfied by response facility $n$ in scenario $\omega$, and $G_{\omega}$ is the Gini coefficient of assignments of relief aid $r$ from response facility $n$ in scenario $\omega$. The optimization model is then

$$(P1') \ \max \sum_{\omega \in [S]} p_{\omega} Q(\mathbf{Y}, \mathbf{P}, \omega) \tag{1}$$

$$\text{s.t.} \sum_{r \in \mathcal{R}} \rho_r P_{rn} \le \sum_{\ell \in \mathcal{L}} \kappa^{\mathrm{resp}}_{\ell} Y_{\ell n} \quad \forall n \in \mathcal{N} \tag{2}$$

$$\sum_{n \in \mathcal{N}} P_{rn} \le \theta^{\max}_{r} \quad \forall r \in \mathcal{R} \tag{3}$$

$$P_{rn} \ge \theta^{\min}_{rn} \sum_{\ell \in \mathcal{L}} Y_{\ell n} \quad \forall r \in \mathcal{R}, n \in \mathcal{N} \tag{4}$$

$$\sum_{\ell \in \mathcal{L}} Y_{\ell n} \leq 1 \quad \forall n \in \mathcal{N} \tag{5}$$

$$\sum_{\substack{r \in \mathcal{R} \\ n \in \mathcal{N}}} c^{\mathrm{p}}_{rn} P_{rn} + \sum_{\substack{\ell \in \mathcal{L} \\ n \in \mathcal{N}}} c^{\mathrm{o}}_{\ell n} Y_{\ell n} \leq \eta \tag{6}$$

$$Y_{\ell n} \in \{0, 1\} \quad \forall \ell \in \mathcal{L}, n \in \mathcal{N} \tag{7}$$

$$P_{rn} \geq 0 \quad \forall r \in \mathcal{R}, n \in \mathcal{N} \tag{8}$$

$$\text{where } Q(\mathbf{Y}, \mathbf{P}, \omega) = \max \sum_{\substack{r \in \mathcal{R} \\ a \in \mathcal{A} \\ n \in \mathcal{N}}} u_{ran\omega} X_{ran\omega} (1 - G_{\omega}) \tag{9}$$

$$\text{s.t.} \sum_{a \in \mathcal{A}} X_{ran\omega} d_{ra\omega} \leq P_{rn} \quad \forall r \in \mathcal{R}, n \in \mathcal{N} \tag{10}$$

$$\sum_{n \in \mathcal{N}} X_{ran\omega} \leq 1 \quad \forall r \in \mathcal{R}, a \in \mathcal{A} \tag{11}$$

$$\sum_{\substack{a \in \mathcal{A} \\ n \in \mathcal{N}}} c^{\mathrm{d}}_{an} \frac{\rho_r}{\kappa^{\mathrm{v}}} d_{ra\omega} X_{ran\omega} \leq \eta' \tag{12}$$

$$X_{ran\omega} \geq 0 \quad \forall r \in \mathcal{R}, a \in \mathcal{A}, n \in \mathcal{N}. \tag{13}$$

The objective function (1) maximizes only the recourse function since there is no first-stage cost. Constraint (2) ensures that relief supplies can be prepositioned at a node $n$ only if a response facility is established at that node and that this prepositioning respects the space limitation of the response facility. Constraint (3) puts an upper bound onthe quantity of relief aid of type $r$ available for prepositioning. This depends largely on supply amounts accumulated through donations and limits on stockpiles for other reasons such as perishability. It is generally agreed a priori between private suppliers and public bodies or non-governmental organisations in charge of disaster relief operations. Constraint (4) maintains that if the decision is made to preposition relief aid of type $r$ at a response facility located at node $n$, a minimum quantity of that relief aid is prepositioned. This prevents operational inefficiencies, such as frequent loading and unloading, associated with having small quantities of supplies stocked at facilities. Constraint (5) restricts each node $n$ to only one type or size of response facility. Constraint (6) defines the pre-disaster financial budgets for carrying out prepositioning activities. Constraints (7) and (8) specify the domains of the first-stage decision variables.

The recourse function given in (9) reflects two crucial concepts in humanitarian logistics, effectiveness and equity. The effectiveness measure is $\sum_{r\in\mathcal{R},a\in\mathcal{A},n\in\mathcal{N}} u_{ran\omega}X_{ran\omega}$, which is the total utility of relief aid assignments, and the equity measure is $1 - G_\omega$, where the Gini coefficient $G_\omega$ is a popular measure of inequity. It is now well-established that equity, whilst being under-studied in humanitarian logistics, is an essential consideration in humanitarian operations (Alem et al. 2016; Sabbaghtorkan et al. 2020; Alem et al. 2022; Çankaya et al. 2019). Given a scenario $s$, the effectiveness function determines the extent to which the established response facilities cover victim needs, whereas the equity function measures the extent to which the prepositioned stock of relief aids is fairly allocated amongst affected areas. The objective function follows the rationale developed in (Eisenhandler and Tzur 2018) with the additional contribution of a novel utility function tailored for the Brazilian case. The explicit form of this function is detailed in the next paragraph. Constraint (10) ensures that assignments of relief aid $r$ from node $n$ to satisfy victim needs do not exceed the available prepositioned supplies at $n$. Constraint (11) guarantees that the total amount of relief aid $r$ allocated to affected area $a$ does not exceed victim needs. Constraint (12) bounds the post-disaster expenses for carrying out the relief aid assignments. Finally, constraints (13) specify the domain of the second-stage decision variables.

### 6.1 Devising an equity-effectiveness trade-off measure

The proposed utility function in (9) is built upon the effectiveness principle using 1) the socioeconomic *vulnerability* of affected areas, 2) the *accessibility*, reflected by the travel times from response facilities to affected areas, 3) the *criticality* of relief aids in alleviating human suffering, and 4) the victim needs, as follows:

$$u_{ran\omega} = \gamma_a \beta_{an} w_r d_{ra\omega},$$

where $\gamma_a$ represents the socioeconomic vulnerability of affected area $a$, $\beta_{ran}$ is the accessibility associated with assigning the response facility at node $n$ to cover victim needs of type $r$ in affected area $a$, $w_r$ is the criticality of relief aid $r$, and $d_{ra\omega}$ represents the victim needs at affected area $a$ for relief aid $r$ in scenario $\omega$. Notice that, for a given coverage level, say $\bar{X}_{ran\omega}$, a scenario-dependent recourse maximization of $\sum_{r\in\mathcal{R},a\in\mathcal{A},n\in\mathcal{N}} u_{ran\omega}\bar{X}_{ran\omega}$ will favour affected areas that present higher utilities.

In this paper, we adopt the poverty measure based on income poverty (FGT) as a proxy for the socioeconomic vulnerability (following developments in Alem et al. (2021)). Poverty is indeed recognized as one of the main drivers that lead to vulnerability to natural hazards in the sense that it narrows coping and resistance strategies, and it causes the loss of diversification, the restriction of entitlements, and the lack of empowerment (de Sherbinin 2008; The World Bank 2017).

The FGT poverty measure is a popular measure for its simplicity and desirable axiomatic properties (Foster 2010). Let $h_a$ be the total population of affected area $a$ and $h_a^{EP}$, $h_a^{VP}$, $h_a^{AP}$ be the number of extremely poor people,very poor people, and almost poor people in $a$, respectively. The income classes represented by extremely poor, very poor, and almost poor people are assumed to have a per capita household income equal to or less than thresholds given by $t^{EP}$, $t^{VP}$, and $t^{AP}$ per month, respectively. The average incomes of these groups are given by $\iota_a^{EP}$, $\iota_a^{VP}$, and $\iota_a^{AP}$, respectively. Let $\iota_0$ be a poverty line or given threshold for income. The FGT poverty measure for an affected area $a$ is calculated as

$$\gamma_a = \frac{1}{h_a}\left[h_a^{EP}\left(\frac{\iota_0 - \iota_a^{EP}}{\iota_0}\right)^2 + h_a^{VP}\left(\frac{\iota_0 - \iota_a^{VP}}{\iota_0}\right)^2 + h_a^{AP}\left(\frac{\iota_0 - \iota_a^{AP}}{\iota_0}\right)^2\right].$$

To quantify the accessibility, we take into account the response time necessary to cover affected area $a$ from node $n$. Let $\tau_{an}$ be the travel time of any relief aid when the response facility at $n$ is assigned to cover affected area $a$, and assume that supplies must ideally arrive at affected areas within a reference time of $\bar{\tau}$ time units. The accessibility is then defined as $\beta_{an} = 1$, if $\tau_{an} \le \bar{\tau}$; $\beta_{an} = 1 - \frac{\tau_{an} - \bar{\tau}}{\bar{\tau}}$, if $\bar{\tau} < \tau_{an} < 2\bar{\tau}$; otherwise, $\beta_{an} = 0$.

To evaluate the importance or criticality of relief aid $r$, let $NP_r$ be the number of people affected by the shortage of relied aid $r$ and $DT_r$ be the maximum deprivation time per person for relief aid $r$. Therefore, $w'_r = \frac{NP_r}{DT_r}$ represents the importance proportional to the number of people and inversely proportional to the deprivation times. Finally, $w_r = \frac{w'_r}{\sum_{\bar{r} \in R} w'_{\bar{r}}}$ to ensure that $\sum_r w_r = 1$.

For the equity measure, we use the relative mean difference proxy for the Gini Coefficient, following the rationale from Mandell (1991), with envy level expressed in terms of our proposed utility, to obtain the following expression:

$$G_\omega = \frac{\sum_{a \in \mathcal{A}} \sum_{a' \in \mathcal{A}: a' > a} \mid W_{a'} \sum_{r \in \mathcal{R}, n \in \mathcal{N}} u_{ran\omega} X_{ran\omega} - W_a \sum_{r \in \mathcal{R}, n \in \mathcal{N}} u_{ra'n\omega} X_{ra'n\omega} \mid}{\sum_{r \in \mathcal{R}, a \in \mathcal{A}, n \in \mathcal{N}} u_{ran\omega} X_{ran\omega}},$$

where $W_a = \frac{\gamma_a}{\sum_{a' \in \mathcal{A}} \gamma_{a'}}$, leading to the following objective function:

$$\max\left\{\sum_{\substack{r \in \mathcal{R}\\ a \in \mathcal{A}\\ n \in \mathcal{N}}} u_{ran\omega} X_{ran\omega} - \sum_{a \in \mathcal{A}} \sum_{a' \in \mathcal{A}\,:\, a' > a} \mid W_{a'} \sum_{\substack{r \in \mathcal{R}\\ n \in \mathcal{N}}} u_{ran\omega} X_{ran\omega} - W_a \sum_{\substack{r \in \mathcal{R}\\ n \in \mathcal{N}}} u_{ra'n\omega} X_{ra'n\omega} \mid \right\},$$

which can be easily linearized to provide a mixed-integer programming formulation.

### 6.2 Model performances

The Brazilian territory is composed of 26 federal units, namely: Acre (AC), Alagoas (AL), Amapa (AP), Amazonas (AM), Bahia (BA), Ceara (CE), Espirito Santo (ES), Goias (GO), Maranhao (MA), Mato Grosso (MT), Mato Grosso do Sul (MS), Minas Gerais (MG), Para (PA), Paraiba (PB), Parana (PR), Pernambuco (PE), Piaui (PI), Rio de Janeiro (RJ), Rio Grande do Norte (RN), Rio Grande Sul (RS), Rondonia (RO), Roraima (RR), Santa Catarina (SC), Sao Paulo (SP), Sergipe (SE) and Tocantins (TO), each of which is represented by a node. Prepositioning is carried out for six types of relief aid: food, water, personal hygiene, cleaning kits, dormitory kit, and mattress. These are typically acquired by the Brazilian government through yearly procurement biddings (ATA 2017).

The victim needs is evaluated as

$$d_{ra\omega} = \left\lceil \frac{\text{length}}{l_r} \right\rceil \times \lceil \text{coverage}_r \times \text{ victims}_{a\omega} \rceil,$$

where, ‘length’ is the number of days in which victims need to be supplied, $l_r$ shows for how long one unit of aid $r$ can supply the victims, $\text{victims}_{a\omega}$ is the number of homeless and displaced people in area $a$ in scenario $\omega$, and $\text{coverage}_r$ is the number of people covered by one unit of relief aid $r$. The characteristics of the different relief aids are summarized in Table 2.

We use the historical data from 2007–2016 from the Integrated Disaster Information System (S2ID 2018) to estimate the victim needs for ten disaster scenarios with equal nominal probabilities of 0.1.

We consider very small, small, medium, large and very large response facilities whose capacities are 1269, 2538, 5076, 11559 and 22087 $m^3$, respectively. The setup cost for a response facility is assumed to be proportional to the construction cost published in monthly reports by the local Unions of Building Construction Industry. Shipping costs are evaluated based on transportation via medium-sized trucks that cover 2.5 km per litre of diesel, at a cost 3 BRL per litre of diesel. The financial budget is taken as ten percent of the minimum total cost required to meet all victim needs. The budgets for both stages are obtained by solving Model ($P1'$) with a cost-minimization objective. The first- and second-stage budgets are BRL 61,413,460 and BRL 181,496 respectively.

**Table 2** Summary of the relief aid characteristics

| Relief aid | Length | Coverage | Volume in $m^3$ | Prep. capacity[a] | Prep. cost |
|---|---|---|---|---|---|
| Food | 30 | 4 | 0.04 | 7846 | 118.8 |
| Water | 1 | 1 | 0.005 | 219583 | 9.9 |
| Hygiene kits | 30 | 4 | 0.04 | 7847 | 109.8 |
| Cleaning kits | 30 | 4 | 0.03 | 7847 | 120.1 |
| Medical products | 30 | 90 | 0.17 | 313690 | 85.20 |
| Mattress | 365 | 1 | 0.011 | 356 | 135 |

[a] ‘Prep. capacity’ refers to the maximum quantity of each aid $r$ that could be acquired

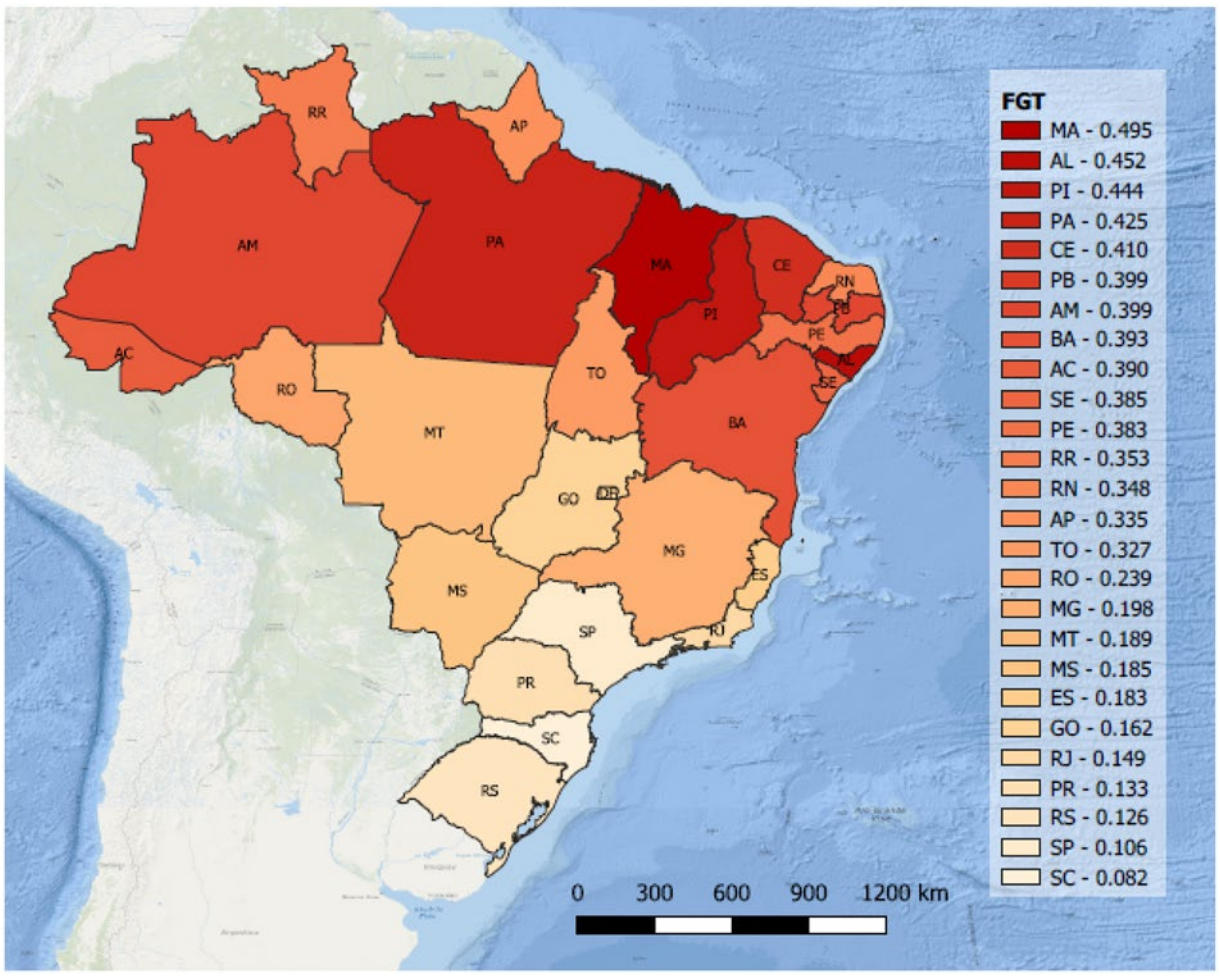

**Fig. 1** FGT index of Brazilian states (Alem et al. 2021)

The socioeconomic data are extracted from the Human Development Atlas published by the United Nations Development Programme at http://www.atlasbrasil.org.br. Figure 1 shows the computed FGT poverty level for all Brazilian states.

Suppose the Kullback–Leibler (K-L) divergence gives an accurate depiction of the decision maker's ambiguity aversion. We use $\Xi = 0.13$ to allow a maximum disaster scenario probability of 0.3. The least-square fits are pictured in Fig. 3.

The LS-PL ($L = U = 5$ is chosen here) mimics excellently the K-L divergence, but provides better practicability. The Moreau-Yosida regularization further improves the fit, which we do not illustrate in Fig. 3 because the difference is visually imperceptible. The SSDs are $4.84 \times 10^{-2}$, $1.48 \times 10^{-4}$, and $1.30 \times 10^{-4}$ for the LS-ICV, LS-PL, and regularized LS-PL, respectively. Table 3 showcases the variations in SSDs of the LS-PL and the regularized LS-PL with the number of pieces. With more pieces, the piecewise linear approximation is closer to the original function. The choice of the number of pieces depends on the order of magnitude of the approximation accuracy that the decision maker wishes to achieve. In our case, we find acceptable an order of magnitude of $10^{-4}$, which is why we choose $L = U = 5$. Figure 3 shows that with this approximation accuracy, the piecewise linear function is almost indistinguishable from the original $f$-divergence function.

**Table 3** Sensitivity analysis on the number of pieces

| | SSD | |
|---|---|---|
| L(=U) | LS-PL | Regularized |
| 1 | $4.72 \times 10^{-2}$ | $5.22 \times 10^{-2}$ |
| 2 | $3.90 \times 10^{-3}$ | $3.90 \times 10^{-3}$ |
| 3 | $8.91 \times 10^{-4}$ | $8.72 \times 10^{-4}$ |
| 4 | $3.16 \times 10^{-4}$ | $3.08 \times 10^{-4}$ |
| 5 | $1.48 \times 10^{-4}$ | $1.30 \times 10^{-4}$ |
| 6 | $7.88 \times 10^{-5}$ | $7.49 \times 10^{-5}$ |
| 7 | $4.76 \times 10^{-5}$ | $4.69 \times 10^{-5}$ |

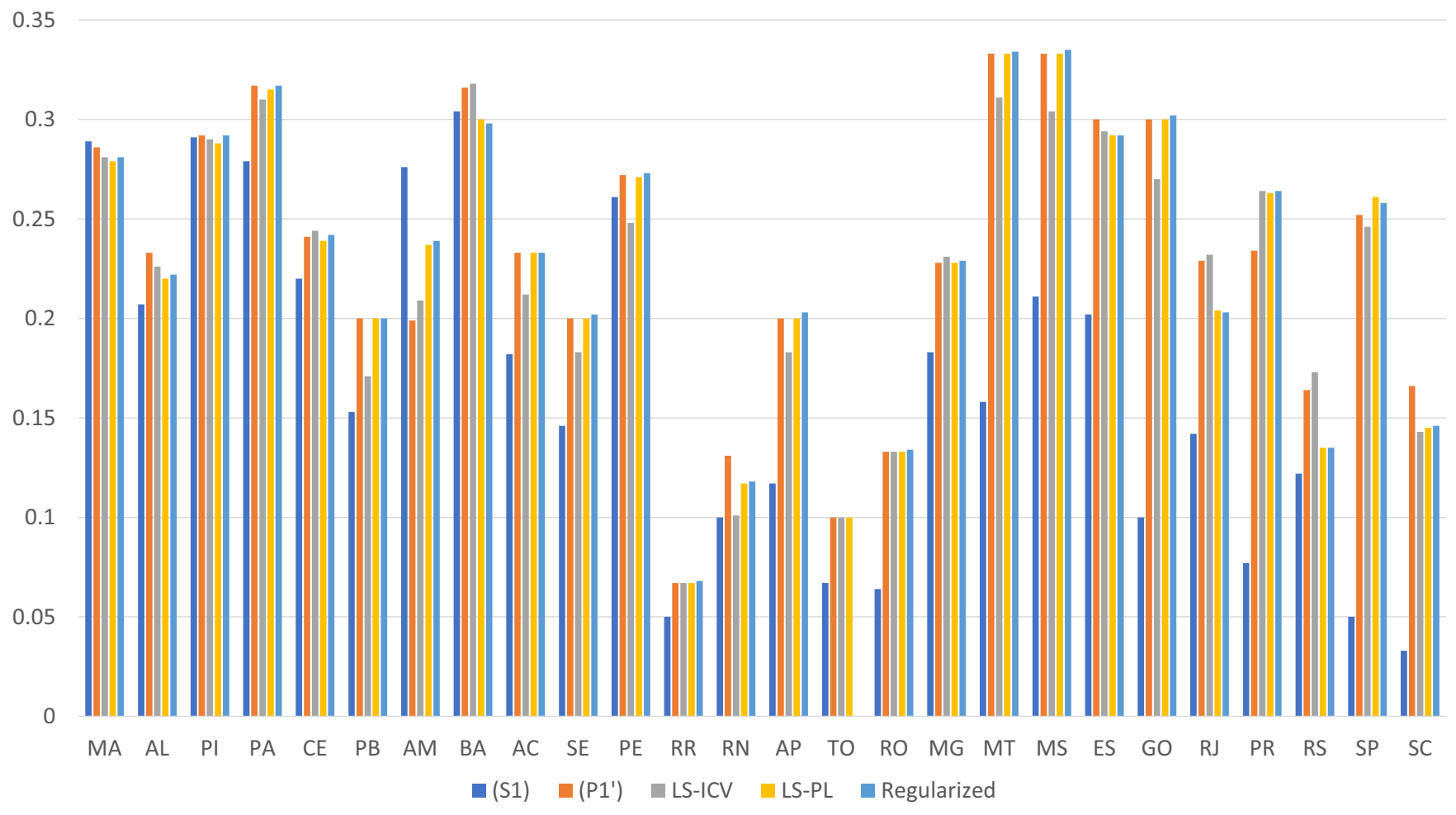

**Fig. 2** Average demand coverage across all federal units for all models

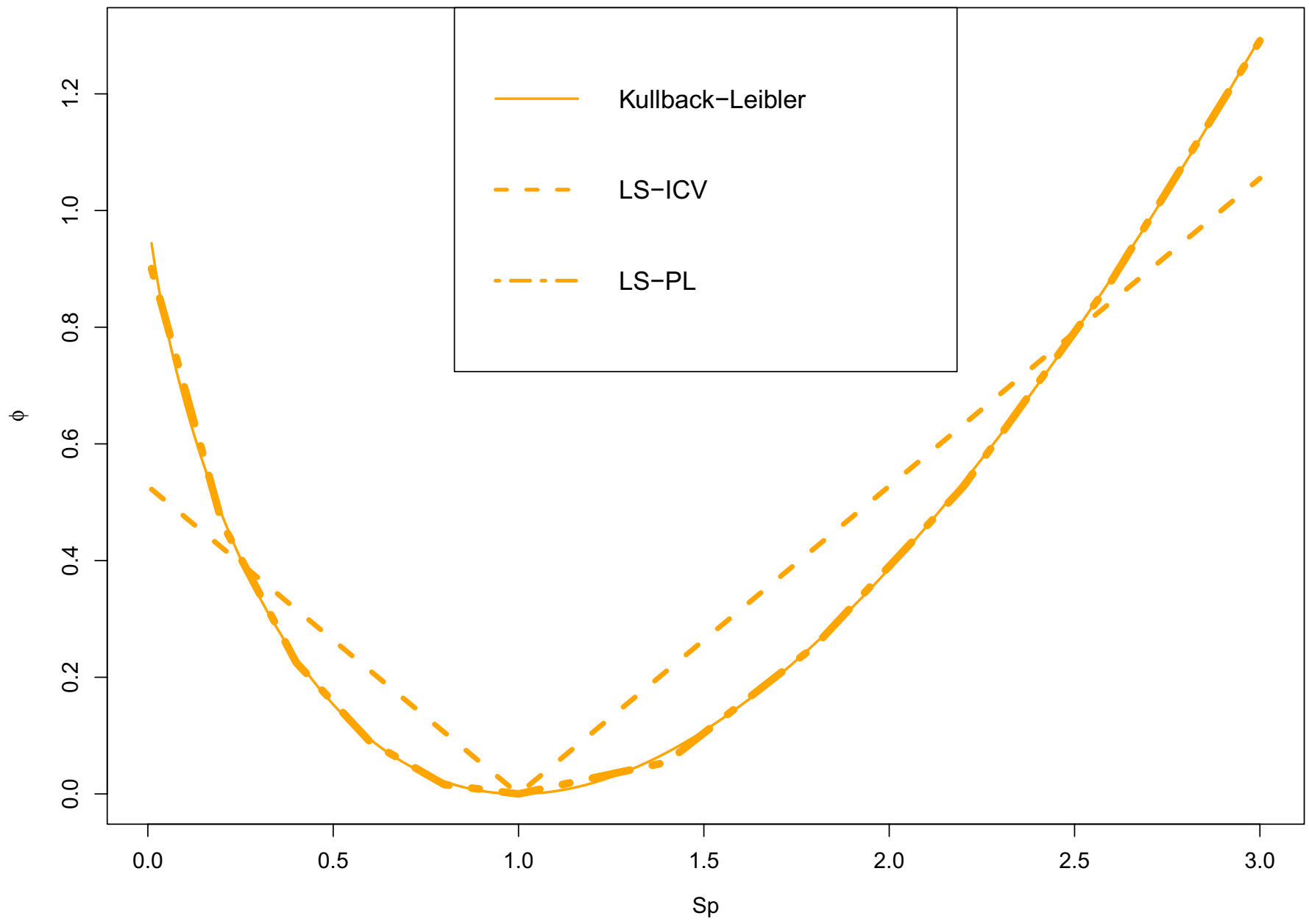

**Fig. 3** Our proposed *f*-divergences

We use CPLEX to solve all models except the model under the K-L divergence, which is mixed-integer with exponential terms. For the latter, we used 3 different solvers in an attempt to get the best solution, BARON, SCIP and DICOPT. The optimal solutions of the various models are illustrated in Table 4. We note that Model ($S1$) is Model ($P1'$) with an effectiveness-only objective ($\sum_{r\in\mathcal{R},a\in\mathcal{A},n\in\mathcal{N}} u_{ran\omega} X_{ran\omega}$).

We see that the model under the K-L divergence fails to converge to the optimal solution after 5 h, whereas the models under our proposed divergences are easily solved. The worst solution time is that of the robust stochastic optimization under the regularized LS-PL, which solves to optimality in under 30 min.

The comparisons of the results of our proposed models with that of a stochastic programming model with an effectiveness-only objective showcase the importance of the Gini coefficient in improving the equity of humanitarian operations. Compared to Model ($S1$), all other models have improved average Gini scores, at the expense of a drop in the effectiveness. We can see the mechanism through which this greater equitability is achieved in Fig. 2.

We see clearly that while the effectiveness-only model ($S1$) provides better average demand coverage for the Amazonas (AM), when equity considerations are factored in, the demand coverage for the Amazonas is sacrificed to obtain better demand coverage in almost every other federal unit. This yields better average demand coverages over all federal units. The average demand coverages are 0.16, 0.23, 0.22, 0.23, 0.22 for Models ($S1$), ($P1'$), LS-ICV, LS-PL, and regularized

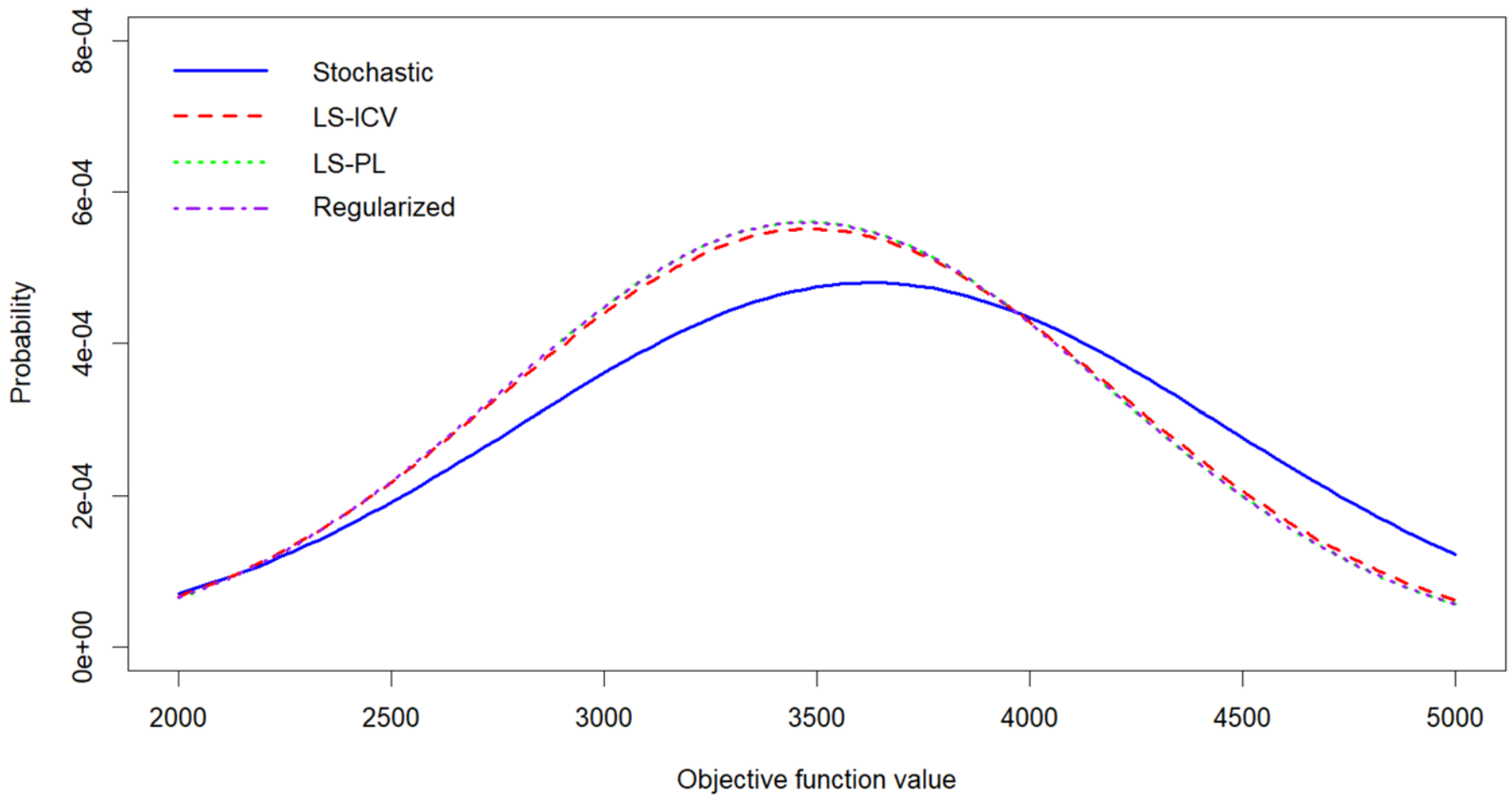

**Fig. 4** Performance under scenario probability variations

**Table 4** Summary of model solutions

| Model | Average demand coverage $\left(\frac{1}{S\|\mathcal{R}\|}\sum_{r\in\mathcal{R},n\in\mathcal{N},\omega\in[S]}X_{ran\omega}\right)$ | | | | | | | | | | | | | | |
|---|---|---|---|---|---|---|---|---|---|---|---|---|---|---|---|
| | AC | AL | AP | AM | BA | CE | ES | GO | MA | MT | MS | MG | PA | PB | PR |
| $(S1)$ | 0.182 | 0.207 | 0.117 | 0.276 | 0.304 | 0.220 | 0.202 | 0.100 | 0.289 | 0.158 | 0.211 | 0.183 | 0.279 | 0.153 | 0.077 |
| $(P1')$ | 0.233 | 0.233 | 0.200 | 0.199 | 0.316 | 0.241 | 0.300 | 0.300 | 0.286 | 0.333 | 0.333 | 0.228 | 0.317 | 0.200 | 0.234 |
| LS-ICV | 0.212 | 0.226 | 0.183 | 0.209 | 0.318 | 0.244 | 0.294 | 0.270 | 0.281 | 0.311 | 0.304 | 0.231 | 0.310 | 0.171 | 0.264 |
| LS-PL | 0.233 | 0.220 | 0.200 | 0.237 | 0.300 | 0.239 | 0.292 | 0.300 | 0.279 | 0.333 | 0.333 | 0.228 | 0.315 | 0.200 | 0.263 |
| Regularized | 0.233 | 0.222 | 0.203 | 0.239 | 0.298 | 0.242 | 0.292 | 0.302 | 0.281 | 0.334 | 0.335 | 0.229 | 0.317 | 0.200 | 0.264 |
| Original $\phi$ | Optimality gap = 12.6% after 5 h running time | | | | | | | | | | | | | | |

| Model | Average demand coverage $\left(\frac{1}{S\|\mathcal{R}\|}\sum_{r\in\mathcal{R},n\in\mathcal{N},\omega\in[S]}X_{ran\omega}\right)$ | | | | | | | | | | | EFF. | E.G. | S.T. |
|---|---|---|---|---|---|---|---|---|---|---|---|---|---|---|
| | PE | PI | RJ | RN | RS | RO | RR | SC | SP | SE | TO | | | |
| $(S1)$ | 0.261 | 0.291 | 0.142 | 0.100 | 0.122 | 0.064 | 0.050 | 0.033 | 0.050 | 0.146 | 0.067 | $1.270\times10^4$ | 0.214 | 3.0 |
| $(P1')$ | 0.272 | 0.292 | 0.229 | 0.131 | 0.164 | 0.133 | 0.067 | 0.166 | 0.252 | 0.200 | 0.100 | $1.064\times10^4$ | 0.282 | 27.0 |
| LS-ICV | 0.248 | 0.290 | 0.232 | 0.101 | 0.173 | 0.133 | 0.067 | 0.143 | 0.246 | 0.183 | 0.100 | $1.073\times10^4$ | 0.267 | 22.0 |
| LS-PL | 0.271 | 0.288 | 0.204 | 0.117 | 0.135 | 0.133 | 0.067 | 0.145 | 0.261 | 0.200 | 0.100 | $1.186\times10^4$ | 0.259 | 34.0 |
| Regularized | 0.273 | 0.292 | 0.203 | 0.118 | 0.135 | 0.134 | 0.068 | 0.146 | 0.258 | 0.202 | 0.100 | $1.189\times10^4$ | 0.259 | 1550.5 |

EFF. is the average effectiveness, defined by $(1/S)\sum_{r\in\mathcal{R},a\in\mathcal{A},n\in\mathcal{N},\omega\in[S]}u_{ran\omega}X_{ran\omega}$

E.G. is the average Gini, defined as $(1/S)\sum_{\omega\in[S]}(1-G_\omega)$

S.T. is the solution time in seconds

Model $(S1)$ is Model $(P1')$ with an effectiveness-only objective $(\sum_{r\in\mathcal{R},a\in\mathcal{A},n\in\mathcal{N}}u_{ran\omega}X_{ran\omega})$

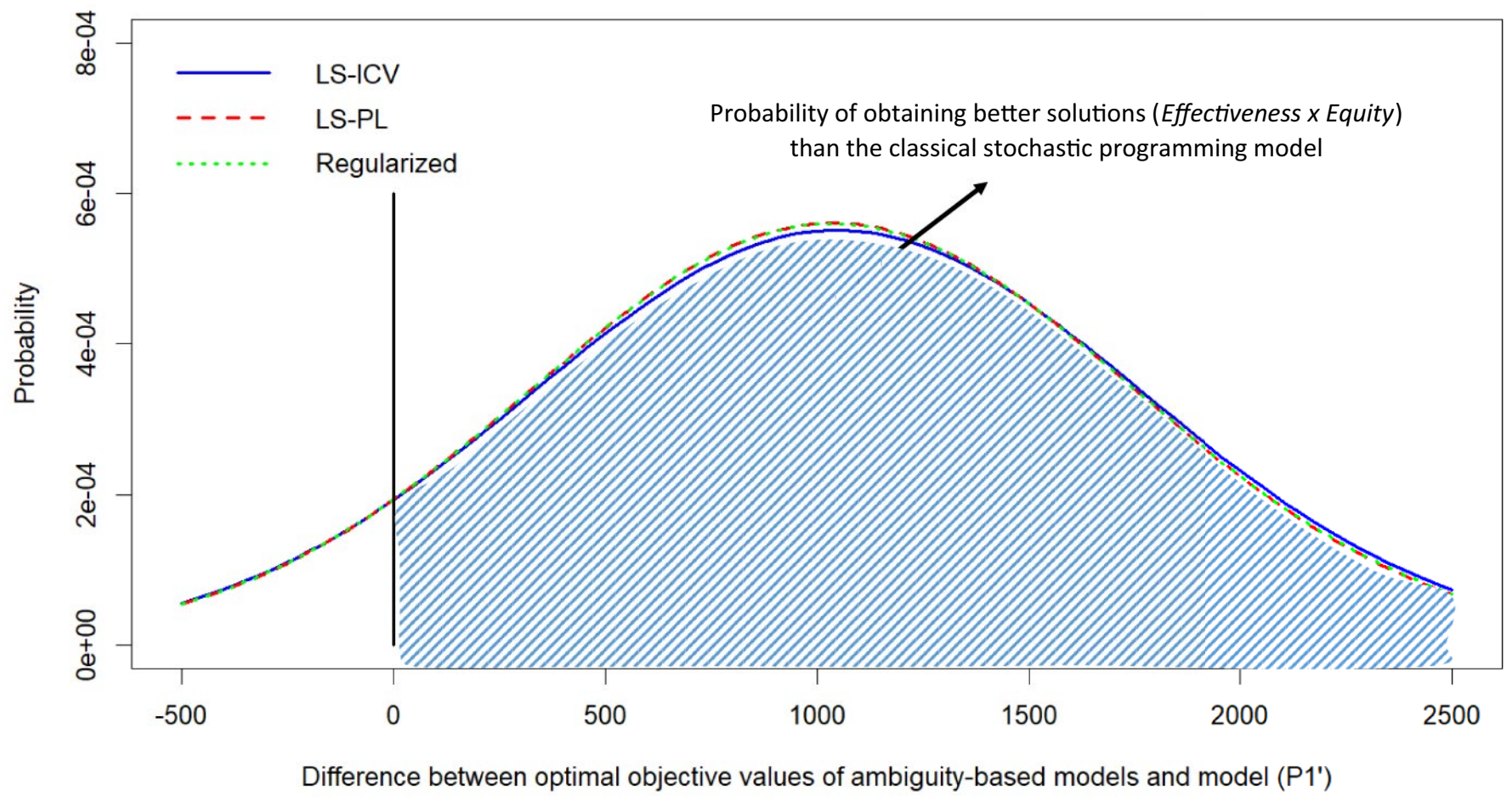

**Fig. 5** Distribution of differences between optimal objective values of ambiguity-based models and model ($P1'$). Positive values indicate that ambiguity-based models have higher (i.e. better) *effectiveness* × *equity* measure in comparison to the classical stochastic programming model

**Table 5** Descriptive statistics of performance improvements of robust stochastic optimization models over the classical stochastic programming model

| | LS-ICV | LS-PL | Regularized |
|---|---|---|---|
| Average | 1048 | 1038 | 1038 |
| Worst-case | −309 | −282 | −279 |
| Best-case | 2671 | 2576 | 2562 |
| Stdev | 731 | 719 | 719 |
| # of scenarios with improvement | 45 | 45 | 45 |
| # of scenarios with deterioration | 5 | 5 | 5 |

LS-PL, respectively. It also yields lower standard deviations of demand coverages, with 0.084, 0.072, 0.071, 0.074, 0.083 for Models ($S1$), ($P1'$), LS-ICV, LS-PL, and regularized LS-PL, respectively. Therefore, resource is redirected from the Amazonas, where victim needs are disproportionately high, to obtain better average and spread of demand coverage across all federal units.

### 6.3 Impact of ambiguity consideration

The previous section studied equity considerations. In this subsection, we investigate the importance of considering ambiguity in humanitarian operations. To do this, we generate 50 random probability vectors such that probabilities do not exceed 0.3 (to keep within the divergence radius). Figure 4 shows the distribution of optimal second-stage objectives when first-stage decisions are fixed with the optimal solutions of the models used in the previous section.

While the performance distributions in Fig. 4 show the overall performance behaviors of the models, they do not portray the *scenario-by-scenario* comparisons of model performances, e.g., in scenario *s*, which model performs best, and by how much? Fig. 5 illustrates how scenario-by-scenario differences between our robust stochastic optimization models and the classical stochastic programming model are then distributed.

The benefit of incorporating ambiguity in humanitarian logistics planning is that it adds greater robustness to wrong probability estimations. In other words, the ambiguity-based approaches hedge against wrong probability predictions, making the solution less sensitive to initial probability assignments, which is where the classical stochastic programming model fails by design. We observe significant improvements in our simulations. Our robust stochastic optimization models yield performance improvements over the classical stochastic programming model in 45 out of the 50 randomly generate probability scenarios. The improvements in the *effectiveness* × *equity* measure are greater than 1000, on average, which corresponds to around 30% of the overall average performance (details are found in Table 5).

## 7 Practicability under random instances

So far, our models have been “trained" on 10 scenarios representing 10 years of Brazilian disaster impacts. Now, we test the practicability of our approaches over larger numbers of scenarios. Additional scenarios are generated randomly with 10% perturbations around existing ones, so as to maintain the correlation of disaster impacts. Values of $\Xi$ are still computed such that maximum probabilities do not exceed three times the probability in the equiprobable setting. Pieces are recomputed for piecewise linear methods (LS-PL and regularized) when the number of scenarios is changed and the solvers and operating system remain unchanged. Our results are shown in Table 6.

The practicability of our robust stochastic optimization models is evident from the results. As we increase the number of scenarios, the solution times of our linear models (LS-ICV and LS-PL) remain comparable to the classical stochastic programming model. Not surprisingly, the robust stochastic optimization model under Moreau-Yosida regularization of piecewise linear divergences has higher solution times for being mixed-integer second-order conic. As mentioned before, this formulation accepts a sacrifice on practicability in order to ensure the differentiability of the divergence function. In the worst-case, the regularized model takes around 5 h to solve. All our robust stochastic optimization models (LS-ICV, LS-PL and Regularized) perform remarkably better than the robust stochastic optimization model under the original *f*-divergence measure. The latter could not close the optimality gap under 500%, even after 10 h of running time. This confirms our original thesis that *f*-divergence functions are largely impractical to implement, especially on mixed-integer models. Our divergence functions approximate *f*-divergences with great accuracy (especially our piecewise linear and regularized piecewise linear functions), while yielding vastly more practicable models.

**Table 6** Practicability over larger number of scenarios

| Solution time (seconds) | | | | | |
|---|---|---|---|---|---|
| # of Scenarios | $(P1')$ | LS-ICV | LS-PL | Regularised[1] | Original[2] $\phi$ |
| 15 | 38 | 21 | 34 | 2293 | 36000 (*gap* > 500%) |
| 20 | 178 | 93 | 44 | 4910 | 36000 (*gap* > 500%) |
| 25 | 214 | 119 | 79 | 7309 | 36000 (*gap* > 500%) |
| 30 | 403 | 152 | 152 | 10158 | 36000 (*gap* > 500%) |
| 35 | 339 | 311 | 508 | 10680 | 36000 (*gap* > 500%) |
| 40 | 451 | 278 | 669 | 10212 | 36000 (*gap* > 500%) |
| 45 | 437 | 399 | 661 | 18820 | 36000 (*gap* > 500%) |
| 50 | 550 | 377 | 717 | 18088 | 36000 (*gap* > 500%) |
| Effectiveness | | | | | |
| 15 | 10202 | 10368 | 11088 | 11043 | 1828 |
| 20 | 10611 | 9981 | 10465 | 11049 | 1810 |
| 25 | 10539 | 9738 | 10465 | 10901 | 1769 |
| 30 | 10537 | 10416 | 10411 | 10872 | 1791 |
| 35 | 10282 | 10384 | 10384 | 11107 | 1701 |
| 40 | 10374 | 10372 | 10459 | 10885 | 1805 |
| 45 | 9789 | 10427 | 10428 | 10703 | 1764 |
| 50 | 10377 | 10218 | 10785 | 10704 | 1697 |
| Inequity | | | | | |
| 15 | 0.286 | 0.276 | 0.269 | 0.270 | 0.435 |
| 20 | 0.282 | 0.287 | 0.277 | 0.273 | 0.440 |
| 25 | 0.281 | 0.287 | 0.271 | 0.273 | 0.446 |
| 30 | 0.283 | 0.285 | 0.285 | 0.276 | 0.442 |
| 35 | 0.284 | 0.282 | 0.282 | 0.270 | 0.447 |
| 40 | 0.286 | 0.286 | 0.284 | 0.277 | 0.441 |
| 45 | 0.288 | 0.283 | 0.283 | 0.278 | 0.448 |
| 50 | 0.286 | 0.281 | 0.279 | 0.281 | 0.443 |

[1] The regularized model is a mixed integer second-order conic problem and it is solved using CPLEX

[2] The original *f*-divergence model is a general non-linear mixed integer optimization problem (with exponential terms). Solution attempts are made using 3 solvers, BARON, SCIP and DICOPT. The suboptimal solution was obtained from DICOPT

In order to make sure that our practicable models are able to maintain their performance improvements over a larger number of scenarios, we ran them on 50 disaster scenarios and 50 randomly generated probability vectors, obtaining the summarized results given in Fig. 6 and Table 7. To offer a greater challenge to our models, we randomly generated 50 *different* disaster scenarios for each probability vector run. Our robust stochastic optimization models improve the classical stochastic programming model on average, with the best improvements coming from our piecewise linear and regularized piecewise linear divergences.

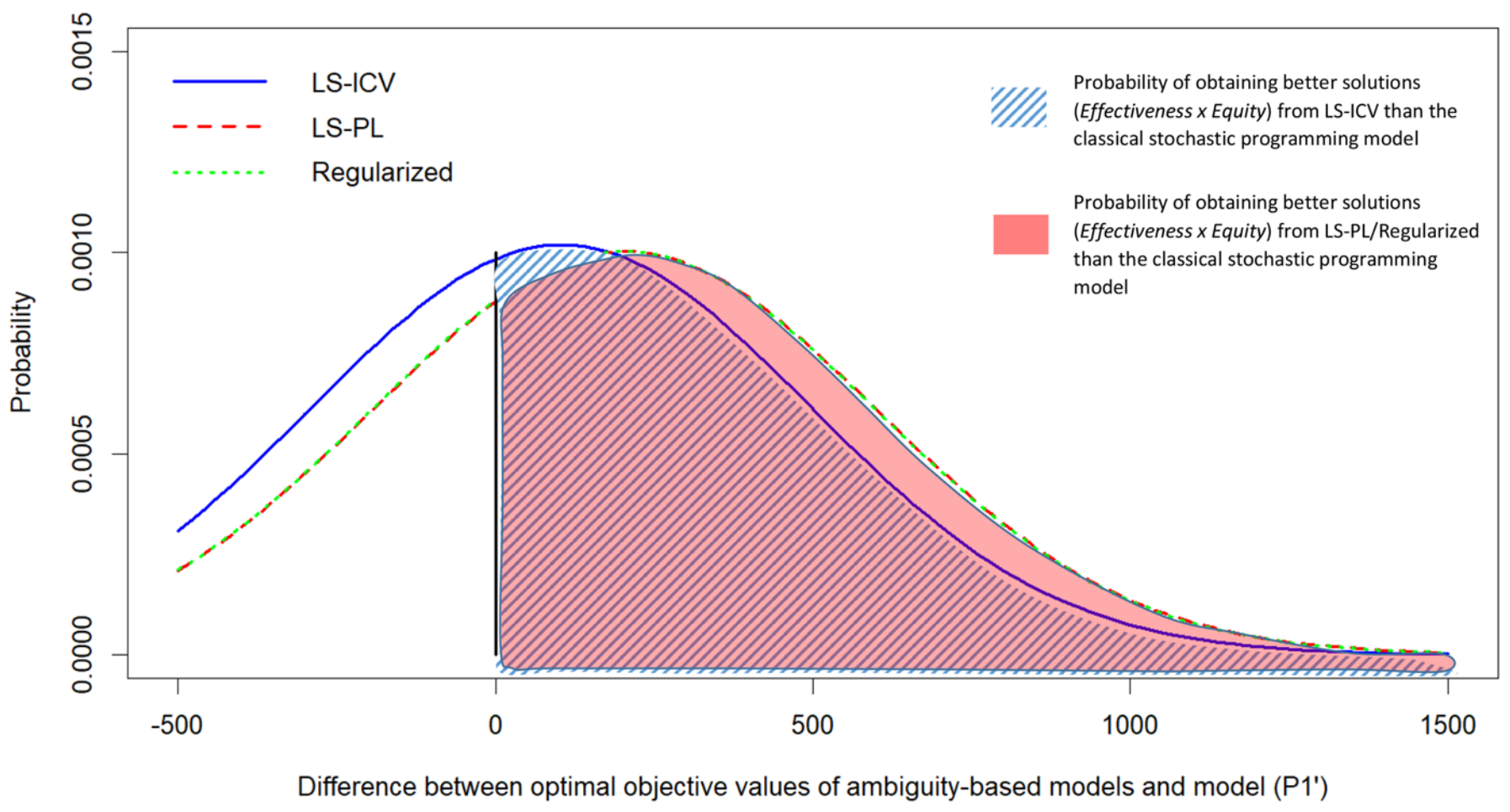

**Fig. 6** Distribution of differences between optimal objective values of ambiguity-based models and model ($P1'$) in the 50-disaster-scenario setting

**Table 7** Descriptive statistics of performance improvements of robust stochastic optimization models over the classical stochastic programming model in the 50-disaster-scenario setting

| | LS-ICV | LS-PL | Regularized |
|---|---|---|---|
| Average | 104 | 204 | 202 |
| Worst-case | −638 | −540 | −542 |
| Best-case | 1344 | 1477 | 1472 |
| Stdev | 395 | 402 | 402 |
| # of scenarios with improvements | 27 | 32 | 32 |
| # of scenarios with deteriorations | 23 | 18 | 18 |

More scenarios experience more improvements than deteriorations and best-case improvements are more pronounced than in the 10-scenario setting.

## 8 Conclusion and future Work

In this paper, we have introduced novel divergence functions that yield practicable robust stochastic optimization counterparts, while providing greater versatility in modeling ambiguity aversions. We have also provided ways in which to tailor these functions such that they mimic existing *f*-divergence functions while offering better practicability. We have applied our approaches to an important problem of disaster preparedness based on a realistic case study of natural hazards in Brazil. Our humanitarian logistics planning model optimizes an effectiveness-equity trade-off, defined with a new utility function and a Gini mean absolute deviation coefficient, both encompassing critical aspects in humanitarian settings, such as equity and vulnerability. With the case study, we showcase 1) the significant improvement in terms of practicability of the robust stochastic optimization counterparts with our proposed divergence functions compared to existing ones, 2) the greater equity of humanitarian response that the objective function enforces and 3) the greater robustness to variations in probability estimations of the resulting plans when ambiguity is considered. In particular, it is remarkable that our findings suggest our robust stochastic optimization models could improve performances on wrong nominal probability estimations in 90% of scenarios, at best, compared to traditional stochastic programming modeling, which is crucial in real-life settings where such improvements mean better allocation of scarce resources to suffering populations. Future research includes enriching our ambiguity sets with CVaR-based metric, in an effort to further reduce the tail of performance deteriorations, thus create even more pronounced performance improvements.

## Appendix: Proofs of theorems, propositions, and lemmas

### Proof of Theorem 1

*Model* $(P1)$ *is equivalent to*

$$\begin{aligned}
\min \quad & \lambda\Xi + \mu^+ - \mu^- + \lambda \sum_{\omega\in[S]} q_\omega \phi^*\left(\frac{\boldsymbol{d}_\omega^\top \boldsymbol{y}_\omega - \mu^+ + \mu^-}{\lambda}\right) \\
\text{s.t.} \quad & \boldsymbol{A}\boldsymbol{x} \geq \boldsymbol{b} \\
& \boldsymbol{D}_\omega \boldsymbol{x} + \boldsymbol{E}\boldsymbol{y}_\omega \geq \boldsymbol{f}_\omega \quad \forall \omega \in [S] \\
& \boldsymbol{x} \in \mathbb{Z}^{N_1'} \times \mathbb{R}^{N_1''} \\
& \boldsymbol{y}_\omega \geq \boldsymbol{0} \quad \forall \omega \in [S] \\
& \lambda, \mu^+, \mu^- \geq 0.
\end{aligned}$$

***Proof*** Proof. Taking the Lagrangian over the ambiguity set and applying the minimax inequality to the model $\sup_{p\in\mathcal{P}} \sum_{\omega\in[S]} p_\omega Q(\boldsymbol{x}, \omega)$, we obtain

$$\begin{aligned}
&\sup_{p\geq 0} \inf_{\lambda,\mu^+,\mu^-\geq 0} \sum_{\omega\in[S]} p_\omega Q(\boldsymbol{x},\omega) + \lambda\left(\Xi - \sum_{\omega\in[S]} q_\omega \phi(\frac{p_\omega}{q_\omega})\right) + (\mu^+ - \mu^-)\left(1 - \sum_{\omega\in[S]} p_\omega\right) \\
&\leq \inf_{\lambda,\mu^+,\mu^-\geq 0} \left\{ \lambda\Xi + \mu^+ - \mu^- + \sup_{p\geq 0} \left\{ \sum_{\omega\in[S]} [p_\omega(Q(\boldsymbol{x},\omega) - \mu^+ + \mu^-) - \lambda q_\omega \phi(\frac{p_\omega}{q_\omega})] \right\} \right\} \\
&= \inf_{\lambda,\mu^+,\mu^-\geq 0} \left\{ \lambda\Xi + \mu^+ - \mu^- + \lambda \sum_{\omega\in[S]} q_\omega \phi^*\left(\frac{Q(\boldsymbol{x},\omega) - \mu^+ + \mu^-}{\lambda}\right) \right\},
\end{aligned}$$

where $\phi^* : \mathbb{R} \to \mathbb{R} \cup \{\infty\}$ is the conjugate of $\phi$. The first inequality follows from the minimax inequality and the second equality follows from the definition of conjugacy. Since $I_{\phi}(\boldsymbol{q},\boldsymbol{q}) = 0$, Slater condition holds and the inequality in the above derivation becomes an equality constraint.

If for any fixed $\omega$, $-\infty < Q(\boldsymbol{x},\omega) < \infty$ for all feasible $\boldsymbol{x}$ and $\phi^*\left(\dfrac{Q(\boldsymbol{x},\omega) - \mu^+ + \mu^-}{\lambda}\right)$ is monotonically increasing in $Q(\boldsymbol{x},\omega)$, combining the first-stage decision, Model $(P1)$ becomes

$$\begin{aligned}
\min \quad & \lambda\Xi + \mu^+ - \mu^- + \lambda \sum_{\omega\in[S]} q_\omega \phi^*\left(\frac{\boldsymbol{d}_\omega^\top \boldsymbol{y}_\omega - \mu^+ + \mu^-}{\lambda}\right)\\
\text{s.t.} \quad & \boldsymbol{A}\boldsymbol{x} \geq \boldsymbol{b}\\
& \boldsymbol{D}_\omega \boldsymbol{x} + \boldsymbol{E}\boldsymbol{y}_\omega \geq \boldsymbol{f}_\omega \quad \forall \omega \in [S]\\
& \boldsymbol{x} \in \mathbb{Z}^{N_1'} \times \mathbb{R}^{N_1''}\\
& \boldsymbol{y}_\omega \geq \boldsymbol{0} \quad \forall \omega \in [S]\\
& \lambda, \mu^+, \mu^- \geq 0.
\end{aligned}$$

□

## Proof of Theorem 2

*For any $\omega \in [S]$, the infimal convolution $\inf_{\sum_{d\in[\mathcal{D}]} r_d = p_\omega/q_\omega}\{\sum_{d\in[\mathcal{D}]} \pi_d \mid \mathcal{D}r_d - 1 \mid\}$, where $\boldsymbol{\pi} > \boldsymbol{0}$, is 1) convex, 2) equal to zero when $p_\omega = q_\omega$ and 3) positive when $p_\omega \neq q_\omega$. Furthermore, when $\pi_d = \frac{1}{\mathcal{D}}$, $\forall d \in [\mathcal{D}]$, this infimal convolution is equivalent to the variation distance $\mid \frac{p_\omega}{q_\omega} - 1 \mid$.*

***Proof*** Proof. The convexity property follows from the fact that the infimal convolution is a convexity-preserving operation. It is clear that, being a positively weighted sum of absolute values, $\inf_{\sum_{d\in[\mathcal{D}]} r_d = p_\omega/q_\omega}\{\sum_{d\in[\mathcal{D}]} \pi_d \mid \mathcal{D}r_d - 1 \mid\} \geq 0$. When $p_\omega = q_\omega$, we can construct a feasible solution with $r_d = \frac{1}{\mathcal{D}}$, $\forall d \in [\mathcal{D}]$, which gives $\sum_{d\in[\mathcal{D}]} \pi_d \mid \mathcal{D}r_d - 1 \mid = 0$, which is the optimal of the convolution. For any $\frac{p_\omega}{q_\omega} \neq 1$, $r_d = \frac{1}{\mathcal{D}}$, $\forall d \in [\mathcal{D}]$, is infeasible and thus, at least one component $\mid \mathcal{D}r_d - 1 \mid$ is positive. Since $\boldsymbol{\pi} > \boldsymbol{0}$, $\inf_{\sum_{d\in[\mathcal{D}]} r_d = p_\omega/q_\omega}\{\sum_{d\in[\mathcal{D}]} \pi_d \mid \mathcal{D}r_d - 1 \mid\} > 0$.

To prove that when $\pi_d = \frac{1}{\mathcal{D}}$, $\forall d \in [\mathcal{D}]$, this infimal convolution is equivalent to the variation distance $\mid \frac{p_\omega}{q_\omega} - 1 \mid$, we rely on the subadditivity of absolute value functions. The infimal convolution

$$\inf_{\sum_{d\in[\mathcal{D}]} r_d = p_\omega/q_\omega} \{\frac{1}{\mathcal{D}} \sum_{d\in[\mathcal{D}]} \mid \mathcal{D}r_d - 1 \mid\}$$
$$= \inf_{r_1,\ldots,r_{[\mathcal{D}-1]}} \{\frac{1}{\mathcal{D}} \sum_{d\in[\mathcal{D}-1]} \mid \mathcal{D}r_d - 1 \mid + \mid \mathcal{D}(\frac{p_\omega}{q_\omega}$$
$$- \sum_{d\in[\mathcal{D}-1]} r_d) - 1 \mid\} \geq \inf_{r_1,\ldots,r_{[\mathcal{D}-1]}} \{\frac{1}{\mathcal{D}} \mid \mathcal{D} \sum_{d\in[\mathcal{D}-1]} r_d - (\mathcal{D}-1) + \mathcal{D}(\frac{p_\omega}{q_\omega}$$
$$- \sum_{d\in[\mathcal{D}-1]} r_d) - 1 \mid\} = \mid \frac{p_\omega}{q_\omega} - 1 \mid \quad .$$

This lower bound is tight and is achieved with the feasible solution $r_d = \frac{p_\omega}{q_\omega \mathcal{D}}$, $\forall d \in [\mathcal{D}]$. □

## Proof of Theorem 3

*Under the divergence* $\inf_{\sum_{d\in[\mathcal{D}]} r_d = p_\omega/q_\omega} \{\sum_{d\in[\mathcal{D}]} \pi_d \mid \mathcal{D}r_d - 1 \mid\}$, *Model* $(P1)$ *is equivalent to the following mixed-integer programming model:*

$$\begin{aligned}
(P3)\ \min \quad & \lambda\Xi + \mu^+ - \mu^- + \sum_{\substack{\omega \in [S] \\ d \in [\mathcal{D}]}} q_\omega z_{\omega d} \\
\text{s.t.} \quad & \boldsymbol{A x} \geq \boldsymbol{b} \\
& \boldsymbol{D}_\omega \boldsymbol{x} + \boldsymbol{E y}_\omega \geq \boldsymbol{f}_\omega \quad \forall \omega \in [S] \\
& z_{\omega d} \geq -\lambda \pi_d \quad \forall \omega \in [S], d \in [\mathcal{D}] \\
& z_{\omega d} \geq \boldsymbol{d}_\omega^\top \boldsymbol{y}_\omega - \mu^+ + \mu^- \quad \forall \omega \in [S], d \in [\mathcal{D}] \\
& \boldsymbol{d}_\omega^\top \boldsymbol{y}_\omega - \mu^+ + \mu^- \leq \lambda \pi_d \quad \forall \omega \in [S], d \in [\mathcal{D}] \\
& \boldsymbol{x} \in \mathbb{Z}^{N_1'} \times \mathbb{R}^{N_1''} \\
& \boldsymbol{y}_\omega \geq \boldsymbol{0}, z_{\omega d} \in \mathbb{R} \quad \forall \omega \in [S], d \in [\mathcal{D}] \\
& \lambda, \mu^+, \mu^- \geq 0 \quad \forall \omega \in [S].
\end{aligned}$$

***Proof*** Proof. The Lagrangian dual of the second stage, $\sup_{\boldsymbol{p}\in\mathcal{P}} \sum_{\omega\in[S]} p_\omega Q(\boldsymbol{x}, \omega)$, of Model $(P1)$ becomes

$$\begin{aligned}
&\inf_{\lambda,\mu^+,\mu^-\geq 0}\left\{\lambda\Xi+\mu^+-\mu^-+\sup_{p\geq 0}\left\{\sum_{\omega\in[S]}[p_\omega(Q(\boldsymbol{x},\omega)-\mu^++\mu^-)\right.\right.\\
&\qquad\left.\left.-\lambda q_\omega \inf_{\sum_{d\in[\mathcal{D}]} r_d=p_\omega/q_\omega}\left\{\sum_{d\in[\mathcal{D}]}\psi_d(r_d)\right\}\right\}\right\}\\
&\quad=\inf_{\lambda,\mu^+,\mu^-\geq 0}\left\{\lambda\Xi+\mu^+-\mu^-+\sup_{p\geq 0}\left\{\sum_{\omega\in[S]}[p_\omega(Q(\boldsymbol{x},\omega)-\mu^++\mu^-)\right.\right.\\
&\qquad\left.\left.-\lambda q_\omega(\psi_1\square\ldots\square\psi_{\mathcal{D}})\left(\frac{p_\omega}{q_\omega}\right)\right\}\right\}\\
&\quad=\inf_{\lambda,\mu^+,\mu^-\geq 0}\left\{\lambda\Xi+\mu^+-\mu^-+\lambda\sum_{\omega\in[S]}q_\omega(\psi_1\square\cdots\square\psi_{\mathcal{D}})^*\left(\frac{Q(\boldsymbol{x},\omega)-\mu^++\mu^-}{\lambda}\right)\right\}\\
&\quad=\inf_{\lambda,\mu^+,\mu^-\geq 0}\left\{\lambda\Xi+\mu^+-\mu^-+\lambda\sum_{\omega\in[S]}q_\omega\sum_{d\in[\mathcal{D}]}\psi_d^*\left(\frac{Q(\boldsymbol{x},\omega)-\mu^++\mu^-}{\lambda}\right)\right\},
\end{aligned} \tag{14}$$

where the first equality results from the definition of infimal convolution, the second equality follows from the definition of conjugacy and the third equality is from the properties of the infimal convolution of proper, convex and lower semicontinuous functions.

We can calculate that

$$\psi_d^*\left(\frac{Q(\boldsymbol{x},\omega)-\mu^++\mu^-}{\lambda}\right)=\begin{cases}\dfrac{Q(\boldsymbol{x},\omega)-\mu^++\mu^-}{D\lambda} & \textit{if } \mid Q(\boldsymbol{x},\omega)-\mu^++\mu^-\mid\leq\lambda D\pi_d\\ -\pi_d \quad \textit{if } Q(\boldsymbol{x},\omega)-\mu^++\mu^-\leq-\lambda D\pi_d.\end{cases}$$

Therefore, $\lambda\psi_d^*(\dfrac{Q(\boldsymbol{x},\omega)-\mu^++\mu^-}{\lambda})$ can be expressed via auxiliary variables in the set $\{z_{\omega d}\in\mathbb{R} : Dz_{\omega d}\geq Q(\boldsymbol{x},\omega)-\mu^++\mu^-, z_{\omega d}\geq-\lambda\pi_d, Q(\boldsymbol{x},\omega)-\mu^++\mu^-\leq\lambda D\pi_d\}$.

Substituting in Model (14), which is both feasible and bounded, and combining with the first-stage model, we obtain Model ($P3$). □

## Proof of Theorem 4

*The minimizer* $\boldsymbol{\pi}^*$ *of* SSD *is such that*

$$\pi_{\underline{d}}^*=\frac{3\Phi}{\mathcal{D}((H-1)^3+1)},$$

*where* $\Phi = \int_0^H \phi(z) \mid z-1 \mid \mathrm{d}z$, *for any one arbitrary* $\underline{d} \in [\mathcal{D}]$ *and* $\pi_d^*$ *is any value greater than* $\pi_{\underline{d}}^*$ *for all* $d \in [\mathcal{D}] \setminus \{\underline{d}\}$.

***Proof*** Proof. WLOG let $\pi_{\underline{d}} = \min_{d\in[\mathcal{D}]}\{\pi_d\}$. Therefore,

$$\int_0^H \left( \inf_{\sum_{d\in[\mathcal{D}]} r_d = z} \left\{ \sum_{d\in[\mathcal{D}]} \pi_d \mid \mathcal{D}r_d - 1 \mid \right\} - \phi(z) \right)^2 \mathrm{d}z$$
$$= \int_0^H \left( \inf_{\{r_d, \forall d \in [\mathcal{D}]\setminus\{\underline{d}\}\}} \left\{ \sum_{d\in[\mathcal{D}]\setminus\{\underline{d}\}} \pi_d \mid \mathcal{D}r_d - 1 \mid + \pi_{\underline{d}} \mid \mathcal{D}\left( z - \sum_{d\in[\mathcal{D}]\setminus\{\underline{d}\}} r_d \right) - 1 \mid \right\} - \phi(z) \right)^2 \mathrm{d}z.$$

It is clear that the optimal value of $r_d$ is $\frac{1}{\mathcal{D}}$, $\forall d \in [\mathcal{D}]\setminus\{\underline{d}\}$. Therefore,

$$\min_{\pi>0} \left\{ \int_0^H \left( \inf_{\{r_d, \forall d \in [\mathcal{D}]\setminus\{\underline{d}\}\}} \left\{ \sum_{d\in[\mathcal{D}]\setminus\{\underline{d}\}} \pi_d \mid \mathcal{D}r_d - 1 \mid + \pi_{\underline{d}} \mid \mathcal{D}\left( z - \sum_{d\in[\mathcal{D}]\setminus\{\underline{d}\}} r_d \right) - 1 \mid \right\} - \phi(z) \right)^2 \mathrm{d}z \right\}$$
$$= \min_{\pi_{\underline{d}}>0} \left\{ \int_0^H \left( \pi_{\underline{d}} \mid \mathcal{D}\left(z - \frac{\mathcal{D}-1}{\mathcal{D}}\right) - 1 \mid - \phi(z) \right)^2 \mathrm{d}z \right\}$$
$$= \min_{\pi_{\underline{d}}>0} \left\{ \int_0^H \left( \pi_{\underline{d}} \mid \mathcal{D}(z-1) \mid - \phi(z) \right)^2 \mathrm{d}z\} = \min_{\pi_{\underline{d}}>0} \{ g(\pi_{\underline{d}}) \right\}.$$

Setting $\frac{\mathrm{d}g}{\mathrm{d}\pi_{\underline{d}}}$ to zero to optimize for $\pi_{\underline{d}}$, we obtain

$$2\pi_{\underline{d}}\mathcal{D}^2 \int_0^H (z-1)^2 \mathrm{d}z - 2\mathcal{D} \int_0^H \phi(z) \mid z-1 \mid \mathrm{d}z = 0$$
$$\Longrightarrow \frac{2\pi_{\underline{d}}\mathcal{D}^2}{3}((H-1)^3 + 1) - 2\mathcal{D}\Phi = 0,$$

where $\Phi = \int_0^H \phi(z) \mid z-1 \mid \mathrm{d}z$. Therefore

$$\pi_{\underline{d}} = \frac{3\Phi}{\mathcal{D}((H-1)^3 + 1)}.$$

Since $r_d = \frac{1}{\mathcal{D}}$, $\forall d \in [\mathcal{D}]\setminus\{\underline{d}\}$ and $\pi_{\underline{d}} = \min_{d\in[\mathcal{D}]}\{\pi_d\}$, $\mid \mathcal{D}r_d - 1 \mid = 0$, $\forall d \in [\mathcal{D}] \setminus \{\underline{d}\}$, it implies that $\pi_d$ can be any value greater than $\pi_{\underline{d}}$, $\forall d \in [\mathcal{D}]\setminus\{\underline{d}\}$.

□

## Proof of Lemma 5

*The LS-ICV of $\phi(z)$ is* $\inf_{\sum_{d\in[\mathcal{D}]} r_d=z}\{\sum_{d\in[\mathcal{D}]}\pi_d \mid \mathcal{D}r_d-1\mid\} = \dfrac{3\Phi}{(H-1)^3+1}\mid z-1\mid$.

***Proof*** Proof. Under the optimal weights, $\pi^*_{\underline{d}} = \dfrac{3\Phi}{\mathcal{D}((H-1)^3+1)}$, for an arbitrary $\underline{d}\in[\mathcal{D}]$ and $\pi^*_d$ being any value greater than $\pi^*_{\underline{d}}$ for all other $d$, the optimal SSD, $\mathrm{SSD}^* = \int_0^H \left(\pi^*_{\underline{d}}\mid \mathcal{D}(z-1)\mid -\phi(z)\right)^2 \mathrm{d}z$, from the proof of Theorem 3. Therefore, $\mathrm{SSD}^* = -\dfrac{3\Phi^2}{(H-1)^3+1}+\Theta$, where $\Theta=\int_0^H \phi^2(z)\mathrm{d}z$. As such, $\mathrm{SSD}^*$ is independent of $\mathcal{D}$ and $\inf_{\sum_{d\in[\mathcal{D}]} r_d=z}\{\sum_{d\in[\mathcal{D}]}\pi^*_d\mid \mathcal{D}r_d-1\mid\} = \dfrac{3\Phi}{(H-1)^3+1}\mid z-1\mid$, which is the infimal convolution when $\mathcal{D}=1$. □

## Proof of Theorem 6

*Under the divergence* $\dfrac{3\Phi}{(H-1)^3+1}\mid \dfrac{p_\omega}{q_\omega}-1\mid, \forall\omega\in[S]$, *Model* $(P1)$ *is equivalent to the following mixed-integer programming model:*

$$
\begin{aligned}
(P4)\ \min\quad & \lambda\Xi+\mu^+-\mu^-+\sum_{\omega\in[S]} q_\omega z_\omega \\
\text{s.t.}\quad & \boldsymbol{Ax}\geq \boldsymbol{b} \\
& \boldsymbol{D}_\omega \boldsymbol{x}+\boldsymbol{E}\boldsymbol{y}_\omega \geq \boldsymbol{f}_\omega \quad \forall\omega\in[S] \\
& z_\omega \geq -\frac{3\lambda\Phi}{(H-1)^3+1} \quad \forall\omega\in[S] \\
& z_\omega \geq \boldsymbol{d}_\omega^{\mathsf{T}}\boldsymbol{y}_\omega-\mu^++\mu^- \quad \forall\omega\in[S] \\
& \boldsymbol{d}_\omega^{\mathsf{T}}\boldsymbol{y}_\omega-\mu^++\mu^- \leq \frac{3\lambda\Phi}{(H-1)^3+1} \quad \forall\omega\in[S] \\
& \boldsymbol{x}\in\mathbb{Z}^{N_1'}\times\mathbb{R}^{N_1''} \\
& \boldsymbol{y}_\omega\geq \boldsymbol{0}, z_\omega\in\mathbb{R} \quad \forall\omega\in[S] \\
& \lambda,\mu^+,\mu^-\geq 0 \quad \forall\omega\in[S].
\end{aligned}
$$

***Proof*** Proof. We begin by deriving the solvable form under our general infimal convolution approximation, $\inf_{\sum_{d\in[\mathcal{D}]} r_d=p_\omega/q_\omega}\{\sum_{d\in[\mathcal{D}]}\psi_d(r_d)\}$. The Lagrangian dual of the second stage, $\sup_{\boldsymbol{p}\in\mathcal{P}}\sum_{\omega\in[S]} p_\omega Q(\boldsymbol{x},\omega)$, of Model $(P1)$ becomes

$$\inf_{\lambda,\mu^+,\mu^-\geq 0}\left\{\lambda\Xi+\mu^+-\mu^-+\sup_{p\geq 0}\left\{\sum_{\omega\in[S]}[p_\omega(Q(\boldsymbol{x},\omega)-\mu^++\mu^-)\right.\right.$$

$$\left.\left.-\lambda q_\omega\inf_{\sum_{d\in[\mathcal{D}]}r_d=p_\omega/q_\omega}\{\sum_{d\in[\mathcal{D}]}\psi_d(r_d)\}\right\}\right\}$$

$$=\inf_{\lambda,\mu^+,\mu^-\geq 0}\left\{\lambda\Xi+\mu^+-\mu^-+\sup_{p\geq 0}\left\{\sum_{\omega\in[S]}[p_\omega(Q(\boldsymbol{x},\omega)-\mu^++\mu^-)\right.\right.$$

$$\left.\left.-\lambda q_\omega(\psi_1\square\cdots\square\psi_{\mathcal{D}})(\frac{p_\omega}{q_\omega})\right\}\right\}$$

$$=\inf_{\lambda,\mu^+,\mu^-\geq 0}\left\{\lambda\Xi+\mu^+-\mu^-+\lambda\sum_{\omega\in[S]}q_\omega(\psi_1\square\cdots\square\psi_{\mathcal{D}})^*(\frac{Q(\boldsymbol{x},\omega)-\mu^++\mu^-}{\lambda})\right\}$$

$$=\inf_{\lambda,\mu^+,\mu^-\geq 0}\left\{\lambda\Xi+\mu^+-\mu^-+\lambda\sum_{\omega\in[S]}q_\omega\sum_{d\in[\mathcal{D}]}\psi_d^*(\frac{Q(\boldsymbol{x},\omega)-\mu^++\mu^-}{\lambda})\right\},$$

where the first equality results from the definition of infimal convolution, the second equality follows from the definition of conjugacy and the third equality is from the properties of the infimal convolution of proper, convex and lower semicontinuous functions. Because of Lemma 1, the above model becomes

$$\inf_{\lambda,\mu^+,\mu^-\geq 0}\left\{\lambda\Xi+\mu^+-\mu^-+\lambda\sum_{\omega\in[S]}q_\omega\chi^*\left(\frac{Q(\boldsymbol{x},\omega)-\mu^++\mu^-}{\lambda}\right)\right\},\tag{15}$$

where $\chi^*$ is the conjugate of $\frac{3\Phi}{(H-1)^3+1}\mid\frac{p_\omega}{q_\omega}-1\mid$.

We can calculate that

$$\chi^*(\frac{Q(\boldsymbol{x},\omega)-\mu^++\mu^-}{\lambda})=\begin{cases}\frac{Q(\boldsymbol{x},\omega)-\mu^++\mu^-}{\lambda} & if\ \mid Q(\boldsymbol{x},\omega)-\mu^++\mu^-\mid\leq\frac{3\lambda\Phi}{(H-1)^3+1}\\ dfrac3\Phi(H-1)^3+1 & if\ Q(\boldsymbol{x},\omega)-\mu^++\mu^-\leq-\frac{3\lambda\Phi}{(H-1)^3+1}.\end{cases}$$

Therefore, $\lambda\chi^*(\frac{Q(\boldsymbol{x},\omega)-\mu^++\mu^-}{\lambda})$ can be expressed via the set of auxiliary variables $\{z_\omega\in\mathbb{R}:z_\omega\geq Q(\boldsymbol{x},\omega)-\mu^++\mu^-,z_\omega\geq-\frac{3\lambda\Phi}{(H-1)^3+1},$ $Q(\boldsymbol{x},\omega)-\mu^++\mu^-\leq\frac{3\lambda\Phi}{(H-1)^3+1}\}$.

Substituting in Model (15), which is both feasible and bounded, and combining with the first-stage model, we obtain Model ($P4$). □

## Proof of Proposition 7

*Suppose a piecewise linear divergence is constructed such that fitting is done separately for ranges $z \le 1$ and $1 \le z \le H$. If L and U pieces with equally spaced intersection points are used for $z \le 1$ and $1 \le z \le H$, respectively, the piecewise linear divergence that fits an existing divergence $\phi(z)$, such that pieces are sequentially least-square fitted starting at $z = 1$, is given by*

$$\mathcal{G}(z) = \begin{cases} \max_{l\in[L]}\left\{\left(3L^3\Psi_l^a - \frac{3L}{2}f_{l+1}^a\left(\frac{l}{L}\right)\right)\left(\frac{l}{L} - z\right) + f_{l+1}^a\left(\frac{l}{L}\right)\right\} & \text{if } 0 \le z \le 1, \\ \max_{u\in[U]}\left\{\left(\frac{3}{\Delta^3}\Psi_u^b - \frac{3}{2\Delta}f_{u-1}^b(1+(u-1)\Delta)(z-1-(u-1)\Delta\right) \right. \\ \left. +f_{u-1}^b(1+(u-1)\Delta)\right\} \\ \text{if } 1 \le z \le H, \end{cases},$$

*where* $f_{L+1}^a(\frac{l}{L}) = 0$, $\Psi_l^a = \int_{(l-1)/L}^{l/L}\phi(z)(\frac{l}{L} - z)\mathrm{d}z$, $\Delta = \frac{H-1}{U}$, $f_0^b(1+(u-1)\Delta) = 0$ *and* $\Psi_u^b = \int_{1+(u-1)\Delta}^{1+u\Delta}\phi(z)(z-1-(u-1)\Delta)\mathrm{d}z$.

***Proof*** Proof. Case $0 \le z \le 1$: The $L^{th}$ piece has the function $a_L(1-z)$ since it must be zero at $z = 1$. Therefore, the SSD is the function $\mathrm{SSD}_L(a_L) = \int_{(L-1)/L}^{1}\left(a_L(1-z) - \phi(z)\right)^2\mathrm{d}z$. Differentiating with respect to $a_L$ and setting to zero to find the minimum SSD, we obtain $\int_{(L-1)/L}^{1}2(1-z)(a_L(1-z)-\phi(z))\mathrm{d}z = 0 \implies a_L = 3L^3\Psi_L$, where $\Psi_L = \int_{(L-1)/L}^{1}\phi(z)(1-z)\mathrm{d}z$.

For $1 \le l \le L-1$, the intersection between $l^{th}$ piece and the $(l+1)^{th}$ piece is at $z = \frac{l}{L}$ since intersection points are equally spaced. The equation for the $l^{th}$ piece is thus $f_l(z) = a_l z - b_l = a_l(\frac{l}{L} - z) + f_{l+1}(\frac{l}{L})$. Therefore, $\mathrm{SSD}_l(a_l) = \int_{(l-1)/L}^{l/L}\left(a_l(\frac{l}{L} - z) + f_{l+1}(\frac{l}{L}) - \phi(z)\right)^2\mathrm{d}z$. Differentiating with respect to $a_l$ and setting to zero, we obtain $\int_{(l-1)/L}^{l/L}2\left(\frac{l}{L} - z\right)\left(a_l\left(\frac{l}{L} - z\right) + f_{l+1}\left(\frac{l}{L}\right) - \phi(z)\right)\mathrm{d}z = 0 \implies \frac{1}{3L^3}a_l + \frac{1}{2L^2}f_{l+1}\left(\frac{l}{L}\right) - \Psi_l = 0$, where $\Psi_l = \int_{(l-1)/L}^{l/L}\phi(z)(\frac{l}{L} - z)\mathrm{d}z$. This implies that $f_l(z) = (3L^3\Psi_l - \frac{3L}{2}f_{l+1}(\frac{l}{L}))(\frac{l}{L} - z) + f_{l+1}(\frac{l}{L})$ and the piecewise linear function is

$$\max_{l\in[L]}\left\{\left(3L^3\Psi_l-\frac{3L}{2}f_{l+1}\left(\frac{l}{L}\right)\right)\left(\frac{l}{L}-z\right)+f_{l+1}\left(\frac{l}{L}\right)\right\},$$

where $f_{L+1}(\frac{l}{L})=0$ and $\Psi_l=\int_{(l-1)/L}^{l/L}\phi(z)(\frac{l}{L}-z)\mathrm{d}z$.

Case $1\le z\le H$: The first piece has the function $a_1(z-1)$ since it must be zero at $z=1$. Letting $\Delta=\frac{H-1}{U}$, the SSD is the function $\mathrm{SSD}_1(a_1)=\int_1^{1+\Delta}\left(a_1(z-1)-\phi(z)\right)^2\mathrm{d}z$. Differentiating with respect to $a_1$ and setting to zero to find the minimum SSD, we obtain $\int_1^{1+\Delta}2(z-1)(a_1(z-1)-\phi(z))\mathrm{d}z=0\implies a_1=\frac{3}{\Delta^3}\Psi_1$, where $\Psi_1=\int_1^{1+\Delta}\phi(z)(z-1)\mathrm{d}z$.

For $2\le u\le U$, the intersection between $u^{th}$ piece and the $(u-1)^{th}$ piece is at $z=1+(u-1)\Delta$ because intersection points are equally spaced. The equation for the $u^{th}$ piece is $f_u(z)=a_uz-b_u=a_u(z-1-(u-1)\Delta)+f_{u-1}(1+(u-1)\Delta)$. Therefore, $\mathrm{SSD}_u(a_u)=\int_{1+(u-1)\Delta}^{1+u\Delta}\left(a_u(z-1-(u-1)\Delta)+f_{u-1}(1+(u-1)\Delta)-\phi(z)\right)^2\mathrm{d}z$. Differentiating with respect to $a_u$ and setting to zero, we obtain $\int_{1+(u-1)\Delta}^{1+u\Delta}2(z-1-(u-1)\Delta)(a_u(z-1-(u-1)\Delta)+f_{u-1}(1+(u-1)\Delta)-\phi(z))\mathrm{d}z=0$ $\implies\frac{\Delta^3}{3}a_u+\frac{\Delta^2}{2}f_{u-1}(1+(u-1)\Delta)-\Psi_u=0$, where $\Psi_u=\int_{1+(u-1)\Delta}^{1+u\Delta}\phi(z)(z-1-(u-1)\Delta)\mathrm{d}z$. This implies that $f_l(z)=(\frac{3}{\Delta^3}\Psi_u-\frac{3}{2\Delta}f_{u-1}(1+(u-1)\Delta)(z-1-(u-1)\Delta)+f_{u-1}(1+(u-1)\Delta)$ and the piecewise linear function is

$$\max_{u\in[U]}\left\{\left(\frac{3}{\Delta^3}\Psi_u-\frac{3}{2\Delta}f_{u-1}(1+(u-1)\Delta)(z-1-(u-1)\Delta)+f_{u-1}(1+(u-1)\Delta\right)\right\},$$

where $\Delta=\frac{H-1}{U}$, $f_0(1+(u-1)\Delta)=0$ and $\Psi_u=\int_{1+(u-1)\Delta}^{1+u\Delta}\phi(z)(z-1-(u-1)\Delta)\mathrm{d}z$. □

## Proof of Corollary 8

*The minimum SSD when the piecewise linear function $\mathcal{G}(z)$ is used to fit an existing divergence $\phi(z)$ is*

$$\text{SSD}^* = \sum_{p\in[L]} \frac{\left(a_p^*\frac{p}{L} - b_p^* - \phi(\frac{p}{L})\right)^3}{3\left(a_p^* - \phi'(\frac{p}{L})\right)} - \sum_{p\in[L]} \frac{\left(a_p^*\frac{p-1}{L} - b_p^* - \phi(\frac{p-1}{L})\right)^3}{3\left(a_p^* - \phi'(\frac{p-1}{L})\right)}$$
$$+ \sum_{p\in[L+U]\setminus[L]} \frac{\left(a_p^*(1+(p-L)\Delta) - b_p^* - \phi(1+(p-L)\Delta)\right)^3}{3\left(a_p^* - \phi'(1+(p-L)\Delta)\right)}$$
$$- \sum_{p\in[L+U]\setminus[L]} \frac{\left(a_p^*(1+(p-L-1)\Delta) - b_p^* - \phi(1+(p-L-1)\Delta)\right)^3}{3\left(a_p^* - \phi'(1+(p-L-1)\Delta)\right)},$$

*where*

$$a_p^* = \begin{cases} -(3L^3\Psi_p^a - \frac{3L}{2}f_{p+1}^a(\frac{p}{L})) & \textit{if } p \in [L] \\ (\frac{3}{\Delta^3}\Psi_{p-L}^b - \frac{3}{2\Delta}f_{p-L-1}^b(1+(p-L-1)\Delta) & \textit{if } p \in [L+U]\setminus[L], \end{cases}$$

*and*

$$b_p^* = \begin{cases} (\frac{3p}{2} - 1)f_{p+1}^a(\frac{p}{L}) - 3pL^2\Psi_p^a & \textit{if } p \in [L] \\ \frac{3+(p-L-1)\Delta}{\Delta^3}\Psi_{p-L}^b - (1 + \frac{3}{2\Delta} + \frac{3(p-L-1)}{2})f_{p-L-1}^b(1+(p-L-1)\Delta) \\ \textit{if } p \in [L+U]\setminus[L] \end{cases}$$

*and* $f_{L+1}^a(\frac{p}{L}) = 0$, $\Psi_p^a = \int_{(p-1)/L}^{p/L} \phi(z)(\frac{p}{L} - z)\mathrm{d}z$, $\Delta = \frac{H-1}{U}$, $f_0^b(1+(p-L-1)\Delta) = 0$, $\Psi_{p-L}^b = \int_{1+(p-L-1)\Delta}^{1+(p-L)\Delta} \phi(z)(z - 1 - (p-L-1)\Delta)\mathrm{d}z$ *and* $\phi'(m)$ *is the derivative of* $\phi$ *evaluated at point* $m$.

***Proof*** Proof. By definition, the minimum SSD is

$$\text{SSD}^* = \int_0^H \left(\mathcal{G}(z) - \phi(z)\right)^2 \mathrm{d}z,$$

where $\mathcal{G}(z) = \max_{p\in[P]}\{a_p^*z - b_p^*\}$. From Proposition 7, we can see that

$$a_p^* = \begin{cases} -\Big(3L^3\Psi_p^a - \frac{3L}{2} f_{p+1}^a\Big(\frac{p}{L}\Big)\Big) & \textit{if } p \in [L] \\ \Big(\frac{3}{\Delta^3}\Psi_{p-L}^b - \frac{3}{2\Delta} f_{p-L-1}^b (1 + (p-L-1)\Delta\Big) & \textit{if } p \in [L+U] \setminus [L], \end{cases}$$

and

$$b_p^* = \begin{cases} \Big(\frac{3p}{2} - 1\Big) f_{p+1}^a\Big(\frac{p}{L}\Big) - 3pL^2\Psi_p^a \quad \textit{if } p \in [L] \\ \frac{3 + (p-L-1)\Delta}{\Delta^3}\Psi_{p-L}^b - \Big(1 + \frac{3}{2\Delta} + \frac{3(p-L-1)}{2}\Big) f_{p-L-1}^b (1 + (p-L-1)\Delta) \\ \textit{if } p \in [L+U] \setminus [L], \end{cases}$$

where $f_{L+1}^a\Big(\frac{p}{L}\Big) = 0$, $\Psi_p^a = \int_{(p-1)/L}^{p/L} \phi(z)(\frac{p}{L} - z)\mathrm{d}z$, $\Delta = \frac{H-1}{U}$, $f_0^b(1 + (p-L-1)\Delta) = 0$, $\Psi_{p-L}^b = \int_{1+(p-L-1)\Delta}^{1+(p-L)\Delta} \phi(z)(z - 1 - (p-L-1)\Delta)\mathrm{d}z$. Given that pieces are equally spaced and the piecewise linear function is convex, we can rewrite SSD* as

$$\begin{aligned} \mathrm{SSD}^* &= \sum_{p\in[L]} \int_{(p-1)/L}^{p/L} \Big(a_p^* z - b_p^* - \phi(z)\Big)^2 \mathrm{d}z + \sum_{p\in[L+U]\setminus[L]} \int_{1+(p-L-1)\Delta}^{1+(p-L)\Delta} \Big(a_p^* z - b_p^* - \phi(z)\Big)^2 \mathrm{d}z \\ &= \sum_{p\in[L]} \frac{\Big(a_p^* z - b_p^* - \phi(z)\Big)^3}{3\Big(a_p^* - \phi'(z)\Big)}\Bigg|_{(p-1)/L}^{p/L} + \sum_{p\in[L+U]\setminus[L]} \frac{\Big(a_p^* z - b_p^* - \phi(z)\Big)^3}{3\Big(a_p^* - \phi'(z)\Big)}\Bigg|_{1+(p-L-1)\Delta}^{1+(p-L)\Delta} \\ &= \sum_{p\in[L]} \frac{\Big(a_p^* \frac{p}{L} - b_p^* - \phi(\frac{p}{L})\Big)^3}{3\Big(a_p^* - \phi'(\frac{p}{L})\Big)} - \sum_{p\in[L]} \frac{\Big(a_p^* \frac{p-1}{L} - b_p^* - \phi(\frac{p-1}{L})\Big)^3}{3\Big(a_p^* - \phi'(\frac{p-1}{L})\Big)} \\ &\quad + \sum_{p\in[L+U]\setminus[L]} \frac{\Big(a_p^*(1 + (p-L)\Delta) - b_p^* - \phi(1 + (p-L)\Delta)\Big)^3}{3\Big(a_p^* - \phi'(1 + (p-L)\Delta)\Big)} \\ &\quad - \sum_{p\in[L+U]\setminus[L]} \frac{\Big(a_p^*(1 + (p-L-1)\Delta) - b_p^* - \phi(1 + (p-L-1)\Delta)\Big)^3}{3\Big(a_p^* - \phi'(1 + (p-L-1)\Delta)\Big)}. \end{aligned}$$

□

## Proof of Theorem 9

*The SSD of the LS-PL is always less than or equal to that of the LS-ICV.*

***Proof*** Proof. A piecewise linear function is defined by $\max_{p\in[P]}\{a_p z - b_p\}$. Since the LS-ICV is the function $\frac{3\Phi}{(H-1)^3+1} \mid z-1 \mid$, it can be expressed as a piecewise linear function where $a_p = b_p = \frac{3\Phi}{(H-1)^3+1}$, $\forall p \in [L]$ and $a_p = b_p = -\frac{3\Phi}{(H-1)^3+1}$, $\forall p \in [L+U] \setminus [L]$. Since the LS-PL is the piecewise linear function that gives minimal SSD, its SSD is always less than or equal to that of LS-ICV. □

## Proof of Theorem 10

*Under the divergence $\mathcal{G}(\frac{p_\omega}{q_\omega})$, Model $(P1)$ is equivalent to the following mixed-integer programming model:*

$$
\begin{aligned}
(P5)\ \min \quad & \lambda\Xi + \mu^+ - \mu^- + \sum_{\omega\in[S]} q_\omega z_\omega \\
\text{s.t.} \quad & \boldsymbol{Ax} \ge \boldsymbol{b} \\
& \boldsymbol{D}_\omega \boldsymbol{x} + \boldsymbol{E}\boldsymbol{y}_\omega \ge \boldsymbol{f}_\omega \quad \forall \omega \in [S] \\
& z_\omega \ge \frac{b_p^* - b_{p+1}^*}{(a_p^* - a_{p+1}^*)}(\boldsymbol{d}_\omega^\top \boldsymbol{y}_\omega - \mu^+ + \mu^- - \lambda a_p^*) + \lambda b_p^* \quad \forall \omega \in [S], p \in [L+U-1] \\
& \boldsymbol{x} \in \mathbb{Z}^{N_1'} \times \mathbb{R}^{N_1''} \\
& \boldsymbol{y}_\omega \ge \boldsymbol{0}, z_\omega \in \mathbb{R} \quad \forall \omega \in [S] \\
& \lambda, \mu^+, \mu^- \ge 0.
\end{aligned}
$$

***Proof*** Proof. The Lagrangian dual of the second stage $\sup_{p\in\mathcal{P}} \sum_{\omega\in[S]} p_\omega Q(\boldsymbol{x},\omega)$ of model $(P2)$ under the piecewise linear fit becomes

$$
\inf_{\lambda,\mu^+,\mu^-\ge 0} \left\{ \lambda\Xi + \mu^+ - \mu^- + \lambda \sum_{\omega\in[S]} q_\omega \mathcal{G}^*\left(\frac{Q(\boldsymbol{x},\omega) - \mu^+ + \mu^-}{\lambda}\right) \right\},
$$

where $\mathcal{G}^*$ is the conjugate of $\mathcal{G}$.

The conjugate of a piecewise linear function $\max_{p\in[L+U]}\{a_p z - b_p\}$ is also a piecewise linear function given by $\sup_z\{zy - \max_{p\in[L+U]}\{a_p z - b_p\}\} = \frac{b_p - b_{p+1}}{a_p - a_{p+1}}(y - a_p) + b_p$, if $a_p \le y \le a_{p+1}$, $\forall p \in [L+U-1]$. Concisely, the conjugate can be written as $\max_{p\in[L+U]}\{\frac{b_p-b_{p+1}}{a_p-a_{p+1}}(y-a_p)+b_p\}$. Replacing the conjugate formulation in the above model and given the latter's feasibility and boundedness, we obtain

$$\begin{aligned}
&\min\ \lambda\Xi+\mu^+-\mu^-+\sum_{\omega\in[S]}q_\omega z_\omega\\
&\text{s.t. } z_\omega\geq\frac{b^*_p-b^*_{p+1}}{(a^*_p-a^*_{p+1})}(Q(\boldsymbol{x},\omega)-\mu^++\mu^--\lambda a^*_p)+\lambda b^*_p \qquad \forall\omega\in[S],p\in[L+U-1]\\
&\lambda,\mu^+,\mu^-\geq 0\\
&z_\omega\in\mathbb{R}\quad\forall\omega\in[S],
\end{aligned}$$

where, from the definition of $\mathcal{G}$, $f^a_{L+1}(\frac{p}{L})=0$, $\Psi^a_p=\int_{(p-1)/L}^{p/L}\phi(z)(\frac{p}{L}-z)\mathrm{d}z$, $\Delta=\frac{H-1}{U}$, $f^b_0(1+(p-L-1)\Delta)=0$, $\Psi^b_{p-L}=\int_{1+(p-L-1)\Delta}^{1+(p-L)\Delta}\phi(z)(z-1-(p-L-1)\Delta)\mathrm{d}z$,

$$a^*_p=\begin{cases}-\left(3L^3\Psi^a_p-\frac{3L}{2}f^a_{p+1}\left(\frac{p}{L}\right)\right) & \textit{if } p\in[L]\\ \left(\frac{3}{\Delta^3}\Psi^b_{p-L}-\frac{3}{2\Delta}f^b_{p-L-1}(1+(p-L-1)\Delta\right) & \textit{if } p\in[L+U]\setminus[L],\end{cases}$$

and

$$b^*_p=\begin{cases}\left(\frac{3p}{2}-1\right)f^a_{p+1}\left(\frac{p}{L}\right)-3pL^2\Psi^a_p \quad \textit{if } p\in[L]\\ \frac{3+(p-L-1)\Delta}{\Delta^3}\Psi^b_{p-L}-\left(1+\frac{3}{2\Delta}+\frac{3(p-L-1)}{2})f^b_{p-L-1}(1+(p-L-1)\Delta\right)\\ \textit{if } p\in[L+U]\setminus[L].\end{cases}$$

Combining with the first-stage model, we obtain Model ($P5$). □

## Proof of Proposition 11

*If $\mathcal{Z}$ points are sampled uniformly from the range $[l'_{p-1}+\epsilon,l'_p-\epsilon]$, $\forall p\in[L+U]$, where $l'_0=0$, $l'_{L+U}=H$ and $l'_p=l_p$ for $p\in\{2,\cdots,L+U-1\}$, the value of $m$ that makes $\mathcal{Y}$ a least-square fit of $\phi$ is obtained from the following quadratic programming problem*

$$m^*=\arg\min_{m>0}\left\{\sum_{p\in[L+U]}\frac{l'_p-l'_{p-1}-2\epsilon}{\mathcal{Z}}\sum_{i\in[\mathcal{Z}]}(a^*_pz_{ip}-b^*_p-\frac{(a^*_p)^2}{2m}-\phi(z_{ip}))^2\right\},$$

*where* $\epsilon\geq\frac{\max_{p\in[L+U]}\{a^*_p\}}{m^*}$.

***Proof*** Proof. If $\mathcal{Z}$ points are sampled uniformly from the range $[l'_{p-1} + \epsilon, l'_p - \epsilon]$, where $\epsilon \geq \dfrac{\max_{p\in[L+U]}\{a^*_p\}}{m}$, $l'_0 = 0$, $l'_{L+U} = H$ and $l'_p = l_p$ for $p \in \{2, \cdots, L + U - 1\}$, the SSD can be calculated as follows:

$$\begin{aligned}\overline{SSD}_p(m) &= \frac{1}{\mathcal{Z}} \sum_{i\in[\mathcal{Z}]} (a^*_p z_{ip} - b^*_p - \frac{(a^*_p)^2}{2m} - \phi(z_{ip}))^2 = \frac{1}{\mathcal{Z}} \sum_{i\in[\mathcal{Z}]} g(z_{ip}) \\ &\approx E\big(g(Z_p)\big) = \frac{1}{l'_p - l'_{p-1} - 2\epsilon} \int_{l'_{p-1}+\epsilon}^{l'_p-\epsilon} g(Z_p)\mathrm{d}Z_p,\end{aligned}$$

where $z_{ip}$ denotes the $i^{th}$ sample point, $Z_p$ is a random variable, $\forall p \in [L + U]$. The last equality is true because of uniform sampling. Therefore,

$$\int_{l'_{p-1}+\epsilon}^{l'_p-\epsilon} g(Z_p)\mathrm{d}Z_p \approx (l'_p - l'_{p-1} - 2\epsilon)\overline{\mathrm{SSD}}_p(m)$$

and the total SSD calculated from sampled points is

$$\begin{aligned}\mathrm{SSD}(m) &= \sum_{p\in[L+U]} \int_{l'_{p-1}+\epsilon}^{l'_p-\epsilon} g(Z_p)\mathrm{d}Z_p \\ &\approx \sum_{p\in[L+U]} \frac{l'_p - l'_{p-1} - 2\epsilon}{\mathcal{Z}} \sum_{i\in[\mathcal{Z}]} (a^*_p z_{ip} - b^*_p - \frac{(a^*_p)^2}{2m} - \phi(z_{ip}))^2.\end{aligned}$$

□

## Proof of Theorem 12

*Under the divergence* $\mathcal{Y}\left(\dfrac{p_\omega}{q_\omega}\right)$, *Model* $(P1)$ *is equivalent to the following mixed-integer second order conic problem:*

$$
\begin{aligned}
(P6)\ \min\quad & \lambda\Xi + \mu^+ - \mu^- + \sum_{\omega\in[S]} q_\omega\Big(z_\omega + \frac{1}{2m}\tau_\omega\Big)\\
\text{s.t.}\quad & \boldsymbol{A}\boldsymbol{x} \ge \boldsymbol{b}\\
& \boldsymbol{D}_\omega\boldsymbol{x} + \boldsymbol{E}\boldsymbol{y}_\omega \ge \boldsymbol{f}_\omega \quad \forall\omega\in[S]\\
& z_\omega \ge \frac{b^*_p - b^*_{p+1}}{(a^*_p - a^*_{p+1})}(\boldsymbol{d}_\omega^\mathsf{T}\boldsymbol{y}_\omega - \mu^+ + \mu^- - \lambda a^*_p) + \lambda b^*_p \quad \forall\omega\in[S], p\in[L+U-1]\\
& \alpha_\omega \ge \boldsymbol{d}_\omega^\mathsf{T}\boldsymbol{y}_\omega - \mu^+ + \mu^- \quad \forall\omega\in[S]\\
& \tau_\omega + \lambda \ge \left\|\begin{pmatrix}\sqrt{2}\alpha_\omega\\ \tau_\omega\\ \lambda\end{pmatrix}\right\|_2 \quad \forall\omega\in[S]\\
& \boldsymbol{x}\in\mathbb{Z}^{N_1'}\times\mathbb{R}^{N_1''}\\
& \boldsymbol{y}_\omega \ge \boldsymbol{0}, \alpha_\omega, \tau_\omega \ge 0, z_\omega\in\mathbb{R} \quad \forall\omega\in[S]\\
& \lambda, \mu^+, \mu^- \ge 0.
\end{aligned}
$$

***Proof*** Proof. With the Moreau-Yosida regularization, the Lagrangian dual of the second stage $\sup_{\boldsymbol{p}\in\mathcal{P}}\sum_{\omega\in[S]} p_\omega Q(\boldsymbol{x},\omega)$ of model $(P1)$ can be written as

$$
\inf_{\lambda,\mu^+,\mu^-\ge 0}\left\{\lambda\Xi + \mu^+ - \mu^- + \lambda\sum_{\omega\in[S]} q_\omega \mathcal{Y}^*\left(\frac{Q(\boldsymbol{x},\omega) - \mu^+ + \mu^-}{\lambda}\right)\right\},
$$

where $\mathcal{Y}^*(\cdot)$ is the conjugate of the Moreau-Yosida regularization. It is known that

$$
\mathcal{Y}^*(\cdot) = \mathcal{G}^*(\cdot) + \frac{1}{2m}\|\cdot\|^2
$$

Replacing the conjugate formulations in the Lagrangian dual and given the latter's feasibility and boundedness, we obtain

$$
\begin{aligned}
&\min\ \lambda\Xi + \mu^+ - \mu^- + \sum_{\omega\in[S]} q_\omega\left(z_\omega + \frac{(Q(\boldsymbol{x},\omega) - \mu^+ + \mu^-)^2}{2m\lambda}\right)\\
&\text{s.t. } z_\omega \ge \frac{b^*_p - b^*_{p+1}}{(a^*_p - a^*_{p+1})}(Q(\boldsymbol{x},\omega) - \mu^+ + \mu^- - \lambda a^*_p) + \lambda b^*_p \qquad \forall\omega\in[S], p\in[L+U-1]\\
&z_\omega\in\mathbb{R} \qquad \forall\omega\in[S]\\
&\lambda, \mu^+, \mu^- \ge 0.
\end{aligned}
$$

This model is equivalent to the second order conic problem

$$\min\ \lambda\Xi + \mu^+ - \mu^- + \sum_{\omega\in[S]} q_\omega\left(z_\omega + \frac{1}{2m}\tau_\omega\right)$$

$$\text{s.t. } z_\omega \geq \frac{b_p^* - b_{p+1}^*}{(a_p^* - a_{p+1}^*)}(Q(\boldsymbol{x},\omega) - \mu^+ + \mu^- - \lambda a_p^*) + \lambda b_p^* \qquad \forall\omega\in[S], p\in[L+U-1]$$

$$\alpha_\omega \geq Q(\boldsymbol{x},\omega) - \mu^+ + \mu^- \qquad \forall\omega\in[S]$$

$$\tau_\omega + \lambda \geq \left\| \begin{pmatrix} \sqrt{2}\alpha_\omega \\ \tau_\omega \\ \lambda \end{pmatrix} \right\|_2 \qquad \forall\omega\in[S]$$

$$\lambda, \mu^+, \mu^- \geq 0$$

$$z_\omega \in \mathbb{R}, \alpha_\omega, \tau_\omega \geq 0 \qquad \forall\omega\in[S].$$

Combining with the first-stage model, we obtain Model ($P6$). □

**Data availability** The generated during and/or analysed during the current study are available in the Mendeley repository, https://data.mendeley.com/datasets/spzgprdjgj/3